\pdfoutput=1
\RequirePackage{ifpdf}
\ifpdf 
\documentclass[pdftex]{sigma}
\else
\documentclass{sigma}
\fi

\numberwithin{equation}{section}

\newtheorem{Theorem}{Theorem}[section]
\newtheorem{Proposition}[Theorem]{Proposition}
\newtheorem{Corollary}[Theorem]{Corollary}

\newtheorem{Lemma}[Theorem]{Lemma}
\theoremstyle{definition}
\newtheorem{Definition}[Theorem]{Definition}

\newtheorem{Remark}[Theorem]{Remark}

\def \D{{\mathcal{D}}}

\def \wt{\operatorname{wt}}

\DeclareMathOperator{\linspan}{span}
\def \Res{\operatorname{Res}}

\def \End{\operatorname{End}}

\def \Hom{\operatorname{Hom}}
\def \mod{\operatorname{mod}}

\def \<{\langle}
\def \>{\rangle}

\begin{document}

\allowdisplaybreaks

\newcommand{\arXivNumber}{2511.01474}

\renewcommand{\PaperNumber}{113}

\FirstPageHeading

\ShortArticleName{Decreasing Filtrations, $C_{2}$-Algebra and Twisted Modules}

\ArticleName{Decreasing Filtrations, $\boldsymbol{C_{2}}$-Algebra\\ and Twisted Modules}

\Author{Shijie CAO~$^{\rm a}$ and Jiancai SUN~$^{\rm ab}$}

\AuthorNameForHeading{S.~Cao and J.~Sun}

\Address{$^{\rm a)}$~Department of Mathematics, Shanghai University, Shanghai 200444, P.R.~China}
\EmailD{\mail{2524709486@qq.com}, \mail{jcsun@shu.edu.cn}}

\Address{$^{\rm b)}$~Newtouch Center for Mathematics, Shanghai University, Shanghai 200444, P.R.~China}

\ArticleDates{Received June 12, 2026, in final form November 03, 2026; Published online November 13, 2026}

\Abstract{We investigate a question posed by Gaberdiel and Gannon (2009) concerning the relationship between $C_{2}$-algebras and twisted modules. To each twisted module $W$ of a vertex algebra $V$, we first associate a decreasing sequence of subspaces \smash{$\bigl\{E_{n}^{T}(W)\bigr\}_{n\in\mathbb{Z}}$} and demonstrate that the associated graded vector space $\mathrm{gr}_{\mathcal{E}}^{T}(W)$ is a twisted module of vertex Poisson algebra $\mathrm{gr}_{\mathcal{E}}^{T}(V)$. We introduce another decreasing sequence of subspaces \smash{$\bigl\{C_{n}^{T}(W)\bigr\}_{n\in\mathbb{Z}_{\geq2}}$} and establish a connection between \smash{$\bigl\{E_{n}^{T}(W)\bigr\}_{n\in\mathbb{Z}}$} and \smash{$\bigl\{C_{n}^{T}(W)\bigr\}_{n\in\mathbb{Z}_{\geq2}}$}. By utilizing the twisted module $\mathrm{gr}_{\mathcal{E}}^{T}(W)$ of vertex Poisson algebra $\mathrm{gr}_{\mathcal{E}}^{T}(V)$, we prove that for any twisted module $W$ of a vertex algebra $V$, $C_{2}$-cofiniteness implies $C_{n}$-cofiniteness for all $n\geq 2$. Furthermore, we employ $\mathrm{gr}_{\mathcal{E}}^{T}(W)$ to study generating subspaces of $\frac{1}{T}\mathbb{N}$-graded twisted modules of lower truncated $\mathbb{Z}$-graded vertex algebras.}


\Classification{17B69}

\section{Introduction}

In \cite{li-abel}, Li introduced and investigated decreasing filtrations for vertex algebras. To any vertex algebra $V$, a canonical decreasing sequence $\mathcal{E}$ of subspaces is associated and it was proven that the associated graded vector space $\mathrm{gr}_{\mathcal{E}}(V)$ naturally forms an $\mathbb{N}$-graded vertex Poisson algebra. The decreasing sequence $\mathcal{E}$ was related to Zhu's decreasing sequence $\mathcal{C}=\{C_{n}\}_{n\geq 2}$. The degree-zero subspace $E_{0}/E_{1}$ of $\mathrm{gr}_{\mathcal{E}}(V)$ was shown to coincide with Zhu's Poisson algebra $V/C_{2}(V)$. It was demonstrated that for any vertex algebra $V$, if $V$ is $C_{2}$-cofinite, then $V$ is $E_{n}$-cofinite and $C_{n+2}$-cofinite for all $n\geq 0$ (see also \cite{bu-span,gn-rq, nt}). Additionally, \cite{li-abel} established that if~$V$ is a~$C_{2}$-cofinite vertex algebra and $W$ is a $C_{2}$-cofinite $V$-module, then $W$ is $C_{n}$-cofinite for all $n\geq 2$. $\mathrm{gr}_{\mathcal{E}}(V)$ was also utilized to study generating subspaces of certain types for lower truncated $\mathbb{Z}$-graded vertex algebras in \cite{li-abel}.

The Li filtration has played a pivotal role in the theory of vertex algebras. It was employed to clarify the equivalence of the two finiteness conditions on vertex algebras in \cite{a-remark} and to study $W$-algebras in \cite{acl-w}. Li filtration for any SUSY vertex algebra was introduced in \cite{y-susy} where it was proven that the associated graded superspace of the filtration has a structure of a SUSY vertex Poisson algebra.

In \cite{gg-zct}, Gaberdiel and Gannon raised a question about clarifying the relation between $C_{2}$-algebra and twisted modules. In this paper, we will investigate this question using decreasing filtrations. We first introduce and study decreasing filtrations for twisted modules of vertex algebras. {Let $V$ be a vertex algebra and let $g$ be an automorphism of $V$ of finite period $T$.} In order to match our results well with the cases of twisted modules, we {first} construct a decreasing sequence $\mathcal{E}^{T}_{V}$ for each vertex algebra $V$, which is a slight generalization of the canonical decreasing sequence $\mathcal{E}$ in \cite{li-abel}, and demonstrate that the associated graded vector space $\mathrm{gr}_{\mathcal{E}}^{T}(V)$ is naturally a vertex Poisson algebra, where for $n\in\mathbb{Z}$, $E_{n}^{T}(V)$ is linearly spanned by the vectors
\[
u_{-1-k_{1}}^{(1)}u_{-1-k_{2}}^{(2)}\cdots u_{-1-k_{r}}^{(r)}v
\]
for $r\geq 1$, $u^{(1)}, u^{(2)}, \dots, u^{(r)}, v\in V$, $k_{1}, k_{2}, \dots,k_{r}\geq 0$ with $k_{1}+k_{2}+\cdots+k_{r}\geq \frac{n}{T}$.

To any $g$-twisted module $W$ of a vertex algebra $V$ we associate a decreasing sequence $\mathcal{E}_{W}^{T}$ of subspaces $E_{n}^{T}(W)$ for $n\in \mathbb{Z}$, which is linearly spanned by the vectors
\[
u^{(1)}_{-1-k_{1}+\frac{r_{1}}{T}}u^{(2)}_{-1-k_{2}+\frac{r_{2}}{T}}\cdots u^{(s)}_{-1-k_{s}+\frac{r_{s}}{T}}w
\]
for $s\geq1$, $u^{(i)}\in V^{r_{i}}$, $0\leq r_{i}\leq T-1$, $1\leq i\leq s$, $w\in W$, $k_{1}, k_{2}, \dots,k_{s}\geq 0$ with $k_{1}+k_{2}+\cdots+k_{s}-\smash{\frac{r_{1}+r_{2}+\cdots +r_{s}}{T}}\geq \smash{\frac{n}{T}}$. We define the notion of twisted modules for vertex Poisson algebras, which is a generalization of the concept of modules of vertex Poisson algebras. {We then} prove that the associated graded vector space $\mathrm{gr}_{\mathcal{E}}^{T}(W)$ is a twisted module {for the} vertex Poisson algebra~$\mathrm{gr}_{\mathcal{E}}^{T}(V)$.

In this paper, we introduce the decreasing sequence \smash{$\mathcal{C}^{T}_{W}=\bigl\{C_{n}^{T}(W)\bigr\}_{n\in\mathbb{Z}_{\geq 2}}$} {associated with the} twisted module~$W$, where for $n\in \mathbb{Z}_{\geq 2}$,
 \[
C^{T}_{n}(W)=\linspan\bigl\{u_{-n+\frac{p}{T}}w \mid u\in V^{p},\, w\in W,\, p=0, 1, \dots, T-1\bigr\}.
\]
We establish a relationship between the sequences $\mathcal{E}^{T}_{W}$ and $\mathcal{C}^{T}_{W}$. Using the twisted module $\mathrm{gr}_{\mathcal{E}}^{T}(W)$ of vertex Poisson algebra $\mathrm{gr}_{\mathcal{E}}^{T}(V)$, we prove that for any twisted module~$W$ of a vertex algebra~$V$, $C_{2}$-cofiniteness implies $C_{n}$-cofiniteness for all $n\geq 2$. Recently, $C_{n}$-cofinite twisted modules for $C_{2}$-cofinite vertex operator algebras with general automorphisms were studied in~\cite{T}.

As shown in this paper, for certain twisted $V$-modules $W$, both sequences $\mathcal{E}^{T}_{W}$ and $\mathcal{C}^{T}_{W}$ are trivial in the sense that $E_{n}^{T}(W)=C_{n+2}^{T}(W)=W$ for all $n\geq 0$. By using the connection between the two decreasing sequences, we prove that if $V=\coprod_{n\geq t}V_{(n)}$ {is} a lower truncated $\mathbb{Z}$-graded vertex algebra for some $t\in\mathbb{Z}$ and $W=\bigoplus_{n\in \frac{1}{T}\mathbb{N}}W_{(n)}$ is a $\frac{1}{T}\mathbb{N}$-graded $g$-twisted $V$-module, then
\begin{gather*}
C^{T}_{n}(W)\subset \coprod_{k\geq n+t-\frac{2T-1}{T}} W_{(k)}\qquad \text{for $n\geq 2$},\\
E^{T}_{m}(W)\subset \coprod_{k\geq n+t-\frac{2T-1}{T}}W_{(k)}\qquad\text{for $m\geq (n-2)T(T+1)2^{n-3}$.}
\end{gather*}
Consequently,
\[
\bigcap_{n\geq0}E^{T}_{n}(W)=\bigcap_{n\geq2}C^{T}_{n}(W)=0.
\]
In this case, both sequences are filtrations. Furthermore, using this result and $\mathrm{gr}_{\mathcal{E}}^{T}(W)$ we show that if there exists a graded subspace $U$ of $V$ such that $V=U+C_{2}(V)$ and a graded subspace~$M$ of a $\frac{1}{T}\mathbb{N}$-graded $g$-twisted $V$-module $W$ such that $W=M+C_{2}^{T}(W)$, then $U$ and $M$ generate~$W$ with a certain spanning property.

This paper is organized as follows: In Section~\ref{section2}, we review the concepts of vertex algebras and twisted modules. In Section~\ref{section3}, we recall the definition of a vertex Poisson algebra and introduce the notion of twisted modules for vertex Lie algebras and vertex Poisson algebras. We then construct a decreasing sequence $\mathcal{E}^{T}_{V}$, and show that the associated graded vector space $\mathrm{gr}_{\mathcal{E}}^{T}(V)$ is a vertex Poisson algebra. In Section~\ref{section4}, we first construct the decreasing filtration $\mathcal{E}^{T}_{W}$ of twisted modules and the associated graded vector space $\mathrm{gr}_{\mathcal{E}}^{T}(W)$, then we prove $\mathrm{gr}_{\mathcal{E}}^{T}(W)$ is a twisted module of vertex Poisson algebra $\mathrm{gr}_{\mathcal{E}}^{T}(V)$. In Section~\ref{section5}, we introduce the decreasing sequence $\mathcal{C}^{T}_{W}$ and establish the relationship between the sequences $\mathcal{E}^{T}_{W}$ and $\mathcal{C}^{T}_{W}$. Finally, in Section~\ref{section6}, by utilizing the twisted module $\mathrm{gr}_{\mathcal{E}}^{T}(W)$ of vertex Poisson algebra $\mathrm{gr}_{\mathcal{E}}^{T}(V)$, we prove that for any twisted module $W$ of a vertex algebra $V$, $C_{2}$-cofiniteness implies $C_{n}$-cofiniteness for all $n\geq 2$. We also study generating subspaces of $\frac{1}{T}\mathbb{N}$-graded $g$-twisted $V$-modules of lower truncated $\mathbb{Z}$-graded vertex algebras.

Throughout this paper, we denote by $\mathbb{Q}$, $\mathbb{Z}$, $\mathbb{Z}_{\geq 2}$ and $\mathbb{N}$ the sets of rational numbers, integers, integers greater than or equal to 2 and nonnegative integers, respectively.

\section{Preliminaries}\label{section2}

In this section, we recall the notions of vertex algebras and twisted modules.
We begin by recalling the notion of vertex algebra (see \cite{b-va, flm,fhl-voa,ll}).

\begin{Definition}\label{def:vertex-algebra}A \emph{vertex algebra} is a vector space $V$, equipped with a linear map
\begin{align*}
Y(\cdot,x)\colon \ & V\rightarrow \Hom (V,V((x)))\subset (\End V)\bigl[\bigl[x,x^{-1}\bigr]\bigr]\\
& v\ \mapsto Y(v,x)=\sum_{n\in \mathbb{Z}}v_{n}x^{-n-1}\qquad (\mbox{where
}v_{n}\in \End V),
\end{align*}
and equipped with a vector ${\bf 1}\in V$, satisfying the conditions
that for $v\in V$,
\[
Y({\bf 1},x)v=v,\qquad Y(v,x){\bf 1}\in V[[x]]\qquad \mbox{and}\qquad
\lim_{x\rightarrow 0}Y(v,x){\bf 1}=v,
\]
and for $u,v,w\in V$, the \emph{Jacobi identity} holds
\begin{align}
&x_{0}^{-1}\delta\left(\frac{x_{1}-x_{2}}{x_{0}}\right)Y(u, x_{1})Y(v, x_{2})w-x_{0}^{-1}\delta\left(\frac{x_{2}-x_{1}}{-x_{0}}\right)Y(v, x_{2})Y(u, x_{1})w\nonumber\\
&\qquad=x_{2}^{-1}\delta\left(\frac{x_{1}-x_{0}}{x_{2}}\right)Y(Y(u, x_{0})v, x_{2})w,\label{jacobi}
\end{align}
where $\delta(x)=\sum_{n\in \mathbb{Z}}x^{n}$. Elementary properties of {the delta function} can be found in~\cite{ll}. All binomial expressions (here and below) are to be expanded in nonnegative integral powers of the second variable, \smash{$(x+y)^{n}=\sum_{k\in \mathbb{N}}\binom{n}{k}x^{n-k}y^{k}$} with \smash{$\binom{n}{k}=\frac{n(n-1)\cdots(n-k+1)}{k!}$} for $n\in \mathbb{Q}$.
\end{Definition}

Taking $\Res_{x_{0}}$ and $\Res_{x_{1}}$ of (\ref{jacobi}) and equating the related coefficients, we have Borcherds' commutator formula and iterate formula:
\begin{gather*}
[u_{m},v_{n}]w=\sum_{i\in \mathbb{N}}\binom{m}{i}(u_{i}v)_{m+n-i}w,\\
(u_{m}v)_{n}w=\sum_{i\in \mathbb{N}}(-1)^{i}\binom{m}{i}\left(u_{m-i}v_{n+i}w-(-1)^{m}v_{m+n-i}u_{i}w\right)
\end{gather*}
for $u,v,w\in V,m,n\in\mathbb{Z}$.

Define a (canonical) linear operator $\mathcal{D}$ on $V$ by
\[
\mathcal{D}(v)=v_{-2}\mathbf{1}\qquad\mathrm{for}\ v\in V.
\]
Then
\[
Y(v,x)\mathbf{1}={\rm e}^{x\mathcal{D}}v\qquad\mathrm{for}\ v\in V.
\]
Furthermore,
\[
[\mathcal{D},v_n]=(\mathcal{D}v)_n=-nv_{n-1}\qquad\mathrm{for}\ v\in V, \ n\in\mathbb{Z}.
\]

Next we recall the definitions of an automorphism of a vertex algebra $V$ and a $g$-twisted $V$-module for a finite order automorphism $g$ of $V$ (see \cite{dlm-twisted,flm}).

\begin{Definition}\label{dfva} An \emph{automorphism} of a vertex algebra $V$ is a linear isomorphism $g$ of $V$ such that $gY(v, x)g^{-1}=Y(g(v), x)$ for any $v\in V$.
\end{Definition}

Let $g$ be a finite order automorphism of $V$ with a period $T$ in the sense that $T$ is a positive integer such that $g^{T}=1$. Then
\[
V=\bigoplus_{r=0}^{T-1}V^{r},
\]
where $V^r=\{v\in V|g(v)={\rm e}^{\frac{2\pi {\rm i}r}{T}}v\}$ for $r\in \mathbb{Z}$. Note that for $r,s\in \mathbb{Z}$, $V^{r}=V^{s}$ if $r\equiv s\ (\mod \ T)$.

\begin{Definition}\label{dfwgtm} Let $V$ be a vertex algebra. A \emph{$g$-twisted $V$-module} is a vector space $M$ equipped with
a~linear map
\begin{align*}
Y_{M}(\cdot, x)\colon \ & V\rightarrow(\End M)\bigl[\bigl[x^{\frac{1}{T}}, x^{-\frac{1}{T}}\bigr]\bigr],\\
& v\mapsto Y_{M}(v, x)=\sum_{n\in \frac{1}{T}\mathbb{Z}}v_{n}x^{-n-1}\qquad (\mbox{where
}v_{n}\in \End M),
\end{align*}
which satisfies the following conditions: For all $u\in V^{r}$, $v\in V$, $w\in M$ with $0\leq r \leq T-1$,
\begin{gather*}
Y_{M}(u, x)=\sum_{n\in \frac{r}{T}+\mathbb{Z}}u_{n}x^{-n-1},\\
u_{n}w=0\qquad \mbox{for $n\in \frac{r}{T}+\mathbb{Z}$ sufficiently large},\\
Y_{M}(\mathbf{1},x)={\rm Id}_{M},
\end{gather*}
and the \emph{twisted Jacobi identity}
\begin{gather}
 x_{0}^{-1}\delta\left(\frac{x_{1}-x_{2}}{x_{0}}\right)Y_{M}(u, x_{1})Y_{M}(v, x_{2})w-x_{0}^{-1}\delta\left(\frac{x_{2}-x_{1}}{-x_{0}}\right)Y_{M}(v, x_{2})Y_{M}(u, x_{1})w\nonumber\\
\qquad=x_{2}^{-1}\left(\frac{x_{1}-x_{0}}{x_{2}}\right)^{-\frac{r}{T}}\delta\left(\frac{x_{1}-x_{0}}{x_{2}}\right)Y_{M}(Y(u, x_{0})v, x_{2})w.\label{twisted-jacobi}
\end{gather}
\end{Definition}

And for $u\in V^{r},w\in M$, $0\leq r\leq T-1$, we also have
\[
Y_{M}(\D u,x)w=\frac{{\rm d}}{{\rm d}x}Y_{M}(u,x)w,
\]
which implies
\[
(\D u)_{-n+\frac{r}{T}}w=\left(n-\frac{r}{T}\right)u_{-n-1+\frac{r}{T}}w
\]
for $n\in \mathbb{Z}$.

Taking $\Res_{x_{0}}$ of (\ref{twisted-jacobi}) we have the commutator relation in twisted case:
\begin{gather}
 [Y_{M}(u,x_{1}),Y_{M}(v,x_{2})]w
 =\Res_{x_{0}}x_{2}^{-1}\!\left(\frac{x_{1}-x_{0}}{x_{2}}\right)^{-\frac{r}{T}} \delta\left(\frac{x_{1}-x_{0}}{x_{2}}\right)Y_{M}(Y(u,x_{0})v,x_{2})w.\!\!\!\label{twisted-commutator}
\end{gather}
For $u\in V^{r}$, $v\in V^{s}$, $0\leq r,s\leq T-1,$ $m,n\in \mathbb{Z}$, equating the coefficients of $x_1^{-m-r/T-1}x_2^{-n-s/T-1}$ in~(\ref{twisted-commutator}), we get the twisted commutator formula (cf.~\cite{kl-gs}):
\begin{gather}
[u_{m+\frac{r}{T}}, v_{n+\frac{s}{T}}]w=\sum_{i\in \mathbb{N}}\binom{m+\frac{r}{T}}{i}(u_{i}v)_{m+n-i+\frac{r+s}{T}}w.\label{twisted-com}
\end{gather}
From (\ref{twisted-commutator}), there is a nonnegative integer $N$ such that
\begin{gather}
(x_{1}-x_{2})^{N}[Y_{M}(u,x_{1}),Y_{M}(v,x_{2})]w=0.\label{twisted-weakcom}
\end{gather}
Then for $m\in\mathbb{Z}$, $n\in\frac{1}{T}\mathbb{Z}$, $u\in V^{r}$, $v\in V$, $w\in M$, $0\leq r\leq T-1$, we have the twisted iterate formula (cf.~\cite{kl-gs}):
\begin{gather}
(u_{m}v)_{n}w=\sum_{i\in \mathbb{N}}\sum_{m\leq j<N}((-1)^{i}\binom{-\frac{r}{T}}{j-m}\binom{j}{i}u_{j-i+\frac{r}{T}}v_{m+n+i-j-\frac{r}{T}}w\nonumber\\
\hphantom{(u_{m}v)_{n}w=}{} -(-1)^{i+j}\binom{-\frac{r}{T}}{j-m}\binom{j}{i}v_{m+n-i-\frac{r}{T}}u_{i+\frac{r}{T}}w).\label{twisted-ite-1}
\end{gather}

And we know that the twisted Jacobi identity is equivalent to the twisted commutator formula~(\ref{twisted-commutator}) and the following associativity relation:
\begin{gather}
(x_{0}+x_{2})^{l+\frac{r}{T}}Y_{M}(u,x_{0}+x_{2})Y_{M}(v,x_{2})w =(x_{2}+x_{0})^{l+\frac{r}{T}}Y_{M}(Y(u,x_{0})v,x_{2})w\label{twisted-asso}
\end{gather}
 for $u\in V^{r}$, $0\leq r\leq T-1$, $w\in M$, where $l$ is a nonnegative integer such that $x^{l+\frac{r}{T}}Y_{M}(u,x)w$ involves only nonnegative integral powers of $x$.

A vertex algebra $V$ equipped with a $\mathbb{Z}$-grading $V=\coprod_{n\in \mathbb{Z}}V_{(n)}$ is called a \emph{$\mathbb{Z}$-graded vertex algebra} if ${\bf 1}\in V_{(0)}$ and for $u\in V_{(k)}$, $k, m, n\in \mathbb{Z}$,
\[
u_{m}V_{(n)}\subset V_{(n+k-m-1)}.
\]
For $u\in V_{(k)}$ with $k\in \mathbb{Z}$, denote $\operatorname{wt}u=k$. We say that a $\mathbb{Z}$-graded vertex algebra $V=\coprod_{n\in \mathbb{Z}}V_{(n)}$ is \emph{lower truncated} if $V_{(n)}=0$ for $n$ sufficiently small. An \emph{$\mathbb{N}$-graded vertex algebra} is defined in the obvious way.

\begin{Definition}\label{dfagtm}
A \emph{$\frac{1}{T}\mathbb{N}$-graded $g$-twisted $V$-module} is a $g$-twisted module with a $\frac{1}{T}\mathbb{N}$-grading
\[
M=\bigoplus_{n\in \frac{1}{T}\mathbb{N}}M_{(n)},
\]
which satisfies the following
\[
v_{m}M_{(n)}\subseteq M_{(n+k-m-1)}
\]
for homogeneous vector $v\in V_{(k)}$, $k\in \mathbb{Z}$, $m\in \frac{1}{T}\mathbb{Z}$, $n\in \frac{1}{T}\mathbb{N}$.
\end{Definition}

For $w\in M_{(n)}$ with $n\in \frac{1}{T}\mathbb{N}$, we denote $\operatorname{wt}w=n$.
We also need the following formula.

\begin{Lemma}\label{tw-ite-2} Let $V$ be a vertex algebra and let $M$ be a $g$-twisted $V$-module. For $u\in V^{r}$, $v\in V$, $w\in M$, $0\leq r\leq T-1$, $m\in \mathbb{Z}$, $n\in \frac{1}{T}\mathbb{Z}$, there exist nonnegative integers $l$ and $k$ such that
\begin{gather}
(u_{m}v)_{n}w=\sum_{i=0}^{k}\sum_{j\in \mathbb{N}}\binom{-l-\frac{r}{T}}{i}\binom{m+i}{j}(-1)^{j}u_{m+l+i-j+\frac{r}{T}}v_{n-l-i+j-\frac{r}{T}}w.\label{twisted-ite-2}
\end{gather}
\end{Lemma}

\begin{proof}
For $u\in V^{r}$, $v\in V$, $w\in M$, $0\leq r\leq T-1$, from (\ref{twisted-asso}) there exists a nonnegative integer~$l$ such that
\[
(x_{0}+x_{2})^{l+\frac{r}{T}}Y_{M}(u,x_{0}+x_{2})Y_{M}(v,x_{2})w=(x_{2}+x_{0})^{l+\frac{r}{T}}Y_{M}(Y(u,x_{0})v,x_{2})w.
\]
For $m\in \mathbb{Z}$, $n\in \frac{1}{T}\mathbb{Z}$, we have
\begin{align*}
 (u_{m}v)_{n}w
& = \Res_{x_{0}}\Res_{x_{2}}x_{0}^{m}x_{2}^{n}Y_{M}(Y(u,x_{0})v,x_{2})w\\
& = \Res_{x_{0}}\Res_{x_{2}}x_{0}^{m}x_{2}^{n}(x_{2}+x_{0})^{-l-\frac{r}{T}}((x_{2}+x_{0})^{l+\frac{r}{T}}Y_{M}(Y(u,x_{0})v,x_{2})w)\\
& = \Res_{x_{0}}\Res_{x_{2}}\sum_{i\in \mathbb{N}}\binom{-l-\frac{r}{T}}{i}x_{0}^{m+i}x_{2}^{n-l-i-\frac{r}{T}}((x_{2}+x_{0})^{l+\frac{r}{T}}Y_{M}(Y(u,x_{0})v,x_{2})w).
\end{align*}
Let $k$ be a nonnegative integer such that $x_{0}^{m+k}Y(u,x_{0})v\in V[[x_{0}]]$. Then
\begin{align*}
 (u_{m}v)_{n}w
& = \Res_{x_{0}}\Res_{x_{2}}\sum_{i=0}^{k}\binom{-l-\frac{r}{T}}{i}x_{0}^{m+i}x_{2}^{n-l-i-\frac{r}{T}}((x_{2}+x_{0})^{l+\frac{r}{T}}Y_{M}(Y(u,x_{0})v,x_{2})w)\\
& = \Res_{x_{0}}\Res_{x_{2}}\sum_{i=0}^{k}\binom{-l-\frac{r}{T}}{i}x_{0}^{m+i} x_{2}^{n-l-i-\frac{r}{T}}((x_{0}+x_{2})^{l+\frac{r}{T}}\\
& \hphantom{=}{}\times Y_{M}(u,x_{0}+x_{2})Y_{M}(v,x_{2})w)\\
& = \Res_{x_{1}}\Res_{x_{2}}\sum_{i=0}^{k}\binom{-l-\frac{r}{T}}{i}(x_{1}-x_{2})^{m+i} x_{2}^{n-l-i-\frac{r}{T}}x_{1}^{l+\frac{r}{T}}Y_{M}(u,x_{1})Y_{M}(v,x_{2})w\\
& = \Res_{x_{1}}\Res_{x_{2}}\sum_{i=0}^{k}\binom{-l-\frac{r}{T}}{i}\sum_{j\in \mathbb{N}}\binom{m+i}{j}(-1)^{j}x_{1}^{m+l+i-j+\frac{r}{T}}x_{2}^{n-l-i+j-\frac{r}{T}}\\
& \hphantom{=}{}\times Y_{M}(u,x_{1})Y_{M}(v,x_{2})w.
\end{align*}
Then we have
\[
(u_{m}v)_{n}w=\sum_{i=0}^{k}\sum_{j\in \mathbb{N}}\binom{-l-\frac{r}{T}}{i}\binom{m+i}{j}(-1)^{j}u_{m+l+i-j+\frac{r}{T}}v_{n-l-i+j-\frac{r}{T}}w
\]
for $u\in V^{r}$, $v\in V$, $w\in M$, $0\leq r\leq T-1$, $m\in \mathbb{Z}$, $n\in \frac{1}{T}\mathbb{Z}$. This completes the proof.
\end{proof}

\section[Decreasing sequence E(V,T) and the vertex Poisson algebra gr(V,T)]{Decreasing sequence $\boldsymbol{\mathcal{E}^{T}_{V}}$\\ and the vertex Poisson algebra $\boldsymbol{\mathrm{gr}_{\mathcal{E}}^{T}(V)}$}\label{section3}

In this section, we first recall the definition of a~vertex Lie algebra and a~vertex Poisson algebra. In order to match our results well with the decreasing sequences of twisted modules, we then construct a decreasing sequence $\mathcal{E}^{T}_{V}$ for each vertex algebra $V$, which is a slight generalization of the canonical decreasing sequence $\mathcal{E}$ in~\cite{li-abel}, and show that the associated graded vector space $\mathrm{gr}_{\mathcal{E}}^{T}(V)$ is naturally a vertex Poisson algebra. Because the proofs are essentially the same as in the decreasing sequence $\mathcal{E}$ case (see the proofs in \cite[Section~2]{li-abel}), we will omit the proofs.

A vertex algebra $V$ is called a \emph{commutative vertex algebra} if
\begin{gather}
[u_{m}, v_{n}]=0\qquad \mbox{for $u, v\in V$, $m,n\in \mathbb{Z}$.}\label{com-va}
\end{gather}
It is well known (see \cite{b-va,fhl-voa}) that (\ref{com-va}) is equivalent to
\[
u_{n}=0\qquad \mbox{for $u\in V$, $n\geq 0$.}
\]
A commutative vertex algebra exactly amounts to a unital commutative associative algebra equipped with a derivation, which is often called a \emph{differential algebra}.

The following definition of the notion of vertex Lie algebra is due to \cite{k-va,P}.

\begin{Definition}\label{def:vertex-lie-algebra}
A \emph{vertex Lie algebra} is a vector space $V$ equipped with a linear operator $D$ and a linear map
\begin{align*}
\begin{split}
Y_{-}(\cdot, x)\colon\ & V\rightarrow\Hom \bigl(V, x^{-1}V\bigl[x^{-1}\bigr]\bigr),\\
& v\mapsto Y_{-}(v, x)=\sum_{n\geq 0}v_{n}x^{-n-1}
\end{split}
\end{align*}
such that for $u, v\in V$, $m,n\in \mathbb{N}$,
\begin{gather*}
(Dv)_{n}=-nv_{n-1},\qquad
u_{n}v=\sum_{i=0}^{\infty} (-1)^{n+i+1}\frac{1}{i!}D^{i}(v_{n+i}u),\\
[u_{m},v_{n}]=\sum_{i=0}^{m}\binom{m}{i}(u_{i}v)_{m+n-i}.
\end{gather*}
\end{Definition}

\begin{Definition}\label{dfvlaau} An \emph{automorphism} of a vertex Lie algebra $V$ is a linear isomorphism $g$ of $V$ such that $gY_{-}(v,x)g^{-1}=Y_{-}(g(v),x)$ for $v\in V$.\end{Definition}

Let $g$ be a finite order automorphism of {a} vertex Lie algebra $V$ with a period $T$ in the sense that $T$ is a positive integer such that $g^{T}=1$. Then
\[
V=\bigoplus_{r=0}^{T-1}V^{r},
\]
where $V^r=\bigl\{v\in V\mid g(v)={\rm e}^{\frac{2\pi {\rm i}r}{T}}v\bigr\}$ for $r\in \mathbb{Z}$.

Next we give the definition of the notion of twisted module of vertex Lie algebra.

\begin{Definition}\label{dftmvla} A \emph{$g$-twisted module} for a vertex Lie algebra $V$ is a vector space $W$ equipped with a linear
map
\begin{align*}
Y_{-}^{W}(\cdot, x)\colon \ & V\rightarrow(\End W)\bigl[\bigl[x^{\frac{1}{T}}, x^{-\frac{1}{T}}\bigr]\bigr],\\
& v\mapsto Y_{-}^{W}(v,x)=\sum_{n\in \frac{1}{T}\mathbb{Z}}v_{n}x^{-n-1}, \qquad \mbox{$v_{n}\in\End W$},
\end{align*}
which satisfies the following conditions: For a linear operator $D$ and $u\in V^{p}$, $v\in V^{q}$, $w\in W$ with $0\leq p,q\leq T-1$, $m,n\in \mathbb{Z}$,
\begin{gather}
Y_{-}^{W}(u,x)=\begin{cases} \displaystyle \sum_{n\in \mathbb{N}}u_{n}x^{-n-1} & \mbox{for $p=0$},\\
\displaystyle\sum_{n\in \mathbb{Z}}u_{n+\frac{p}{T}}x^{-n-\frac{p}{T}-1} & \mbox{for $1\leq p\leq T-1$},
\end{cases}\nonumber
\\
u_{n+\frac{p}{T}}w=0\qquad \mbox{for $n$ sufficiently large},\nonumber
\\
(Du)_{n+\frac{p}{T}}=\left(-n-\frac{p}{T}\right)u_{n-1+\frac{p}{T}},\label{dftmvla-1}\\
\bigl[u_{m+\frac{p}{T}},v_{n+\frac{q}{T}}\bigr]=\sum_{i\in \mathbb{N}}\binom{m+\frac{p}{T}}{i}(u_{i}v)_{m+n-i+\frac{p+q}{T}}\label{dftmvla-2}.
\end{gather}
\end{Definition}

Note that from Definition~\ref{dftmvla}, we have
\[
u_{-j+\frac{p}{T}}=\frac{1}{\bigl(j-1-\frac{p}{T}\bigr)\bigl(j-2-\frac{p}{T}\bigr)\cdots\bigl(-\frac{p}{T}\bigr)}\bigl(D^{j} u\bigr)_{\frac{p}{T}}
\]
 for $u\in V^{p}$, $j\geq 1$, $1\leq p\leq T-1$.

Recall the following notion of vertex Poisson algebra from \cite{dlm-vpa,fb-va,li-vpa}.

\begin{Definition}\label{dfvp}
A \emph{vertex Poisson algebra} is a commutative vertex algebra $A$, or equivalently,
a~(unital) commutative associative algebra with a derivation $\partial$, equipped with a vertex Lie algebra
structure $(Y_-,\partial)$ such that
\[
Y_-(a,x)\in x^{-1}(\operatorname{Der}A)\bigl[x^{-1}\bigr]\qquad \mbox{for $a\in A$}.
\]\end{Definition}

Next we give the following definition of the notion of twisted module for vertex Poisson algebra.

\begin{Definition}\label{tmvpam}
Let $A$ be a vertex Poisson algebra and let $g$ be an automorphism of $A$
of finite period $T$. A $g$-twisted $A$-module is a vector space $W$ equipped
with a $g$-twisted module structure for the underlying vertex Lie algebra
$(A,Y_-,\partial)$ in the sense of Definition~\ref{dftmvla}, and a bilinear map
\[
{A\otimes W\longrightarrow W,\qquad a\otimes w\longmapsto a\cdot w}
\]
satisfying $\mathbf{1}\cdot w=w$ for every $w\in W$ and the compatibility condition
\begin{gather}
Y^{W}_{-}(u,x)(v\cdot w)=(Y_-(u,x)v)\cdot w+v\cdot Y^{W}_{-}(u,x)w\label{comp}
\end{gather}
for $u,v\in A$ and $w\in W$.\end{Definition}

Below for each vertex algebra, we construct a decreasing sequence \smash{$\mathcal{E}^{T}_{V} = \bigl\{E_{n}^{T}(V)\bigr\}_{n\in \mathbb{Z}}$}, as a~slight generalization of the canonical decreasing sequence $\mathcal{E}$ in~\cite{li-abel}.

\begin{Definition}\label{dftf}Let $V$ be a vertex algebra and let $g$ be an automorphism of $V$ with a period~$T$. Define a sequence \smash{$\mathcal{E}_{V}^{T} =\bigl\{ E_{n}^{T}(V)\bigr\}_{n\in\mathbb{Z}}$} of subspaces of~$V$, where for $n\in\mathbb{Z}$, $E_{n}^{T}(V)$ is linearly spanned by the vectors
\[
u_{-1-k_{1}}^{(1)}u_{-1-k_{2}}^{(2)}\cdots u_{-1-k_{r}}^{(r)}v
\]
for $r\geq 1$, $u^{(1)}, u^{(2)}, \dots, u^{(r)}, v\in V$, $k_{1}, k_{2}, \dots,k_{r}\geq 0$ with $k_{1}+k_{2}+\cdots+k_{r}\geq \frac{n}{T}$.
\end{Definition}

The following are some immediate consequences.

\begin{Lemma}\label{property-en}Let $V$ be a vertex algebra and let $g$ be an automorphism of $V$ with a period $T$. We have
\begin{gather*}
E_{n}^{T}(V)\supset E_{n+1}^{T}(V)\qquad \mbox{for $n\in\mathbb{Z}$,}\\
E_{1+kT}^{T}(V)=E_{2+kT}^{T}(V)=\cdots=E_{T+kT}^{T}(V)\qquad \mbox{for $k\in\mathbb{Z}$,}\\
E_{n}^{T}(V)=V\qquad \mbox{for $n\leq 0$,}\\
u_{-1-k}E_{n}^{T}(V)\subset E_{n+kT}^{T}(V)\qquad \mbox{for $u\in V$, $k\geq 0$, $n\in\mathbb{Z}$.}
\end{gather*}
\end{Lemma}

The following gives a stronger spanning property for $E_{n}^{T}(V)$.

\begin{Lemma}\label{span-en}Let $V$ be a vertex algebra and let $g$ be an automorphism of $V$ with a period $T$. For $n\geq 1$, we have
\begin{gather}
E_{n}^{T}(V)=\linspan\bigl\{u_{-1-i}v\mid u\in V,\, i\geq 1,\, v\in E_{n-iT}^{T}(V)\bigr\}.\label{spann-i}
\end{gather}
Furthermore, for $n\geq 1$, $E_{n}^{T}(V)$ is linearly spanned by the vectors
\[
u_{-1-k_{1}}^{(1)}u_{-1-k_{2}}^{(2)}\cdots u_{-1-k_{r}}^{(r)}v
\]
for $r\geq 1$, $u^{(1)}, u^{(2)}, \dots, u^{(r)}, v\in V$, $k_{1}, k_{2}, \dots, k_{r}\geq 1$ with $k_{1}+k_{2}+\cdots+k_{r}\geq \frac{n}{T}$.
\end{Lemma}

{First,} we have the following special case for $\mathcal{E}_{V}^{T}$.

\begin{Lemma}\label{am-en}Let $V$ be a vertex algebra and let $g$ be an automorphism of $V$ with a period $T$. For $a\in V$, $m,n\in\mathbb{Z}$, we have
\[
a_{m}E_{n}^{T}(V)\subset E_{n-(m+1)T}^{T}(V).
\]
Furthermore,
\[
a_{m}E_{n}^{T}(V)\subset E_{n-(m+1)T+1}^{T}(V)\qquad \mbox{for $m\geq 0$}.
\]
\end{Lemma}

Now we have the following general case.

\begin{Proposition}\label{unv-en} Let $V$ be a vertex algebra and let $g$ be an automorphism of~$V$ with a period~$T$. Let $u\in E_{r}^{T}(V)$, $v\in E_{s}^{T}(V)$ with $r,s\in\mathbb{Z}$. Then
\[
u_{n}v\in E_{r+s-(n+1)T}^{T}(V)\qquad \mbox{for $n\in\mathbb{Z}$}.
\]
Furthermore, we have
\begin{gather}
u_{n}v\in E_{r+s-(n+1)T+1}^{T}(V)\qquad \mbox{for $n\geq 0$}.\label{general-rs0}
\end{gather}
\end{Proposition}

From Proposition \ref{unv-en}, we immediately have the following.

\begin{Theorem}\label{th-vpa} Let $V$ be a vertex algebra, let $g$ be an automorphism of $V$ with a period $T$ and let $\mathcal{E}_{V}^{T}=\{E_{n}^{T}(V)\}_{n\in\mathbb{Z}}$ be the decreasing sequence defined in Definition~{\rm \ref{dftf}} for $V$. Set
\[
\mathrm{gr}^{T}_{\mathcal{E}}(V)=\coprod_{n\in \mathbb{N}}E_{nT}^{T}(V)/E_{(n+1)T}^{T}(V).
\]
Then $\mathrm{gr}^{T}_{\mathcal{E}}(V)$ is a vertex Poisson algebra if we define
\begin{gather}
\bigl(u+E_{(r+1)T}^{T}(V)\bigr)\cdot\bigl(v+E_{(s+1)T}^{T}(V)\bigr)=u_{-1}v+E_{(r+s+1)T}^{T}(V),\label{grv-1}\\
\partial\bigl(u+E_{(r+1)T}^{T}(V)\bigr)=\D(u)+E_{(r+2)T}^{T}(V),\label{grv-2}\\
Y_{-}\bigl(u+E_{(r+1)T}^{T}(V),x\bigr)\bigl(v+E_{(s+1)T}^{T}(V)\bigr)
 =\sum_{n\in \mathbb{N}}\bigl(u_{n}v+E_{(r+s-n+1)T}^{T}(V)\bigr)x^{-n-1}\nonumber
\end{gather}
for $u\in E_{rT}^{T}(V)$, $v\in E_{sT}^{T}(V)$ with $r, s\in \mathbb{N}$.\end{Theorem}

\begin{Remark}\label{avagr}
Note that vertex Poisson algebra $\mathrm{gr}^{T}_{\mathcal{E}}(V)$ in Theorem~\ref{th-vpa} is isomorphic to vertex Poisson algebra $\mathrm{gr}_{\mathcal{E}}(V)$ in~\cite[Theorem 2.12]{li-abel}.
\end{Remark}

Next we recall the sequence $\mathcal{C}$ introduced by Zhu and we then give a relation between Zhu's Poisson algebra and vertex Poisson algebra $\mathrm{gr}^{T}_{\mathcal{E}}(V)$.
\begin{Definition}\label{dfCn}
Let $V$ be a vertex algebra. For $n\geq 2$, define
\[
C_{n}(V)=\linspan\{u_{-n}v\mid u,v\in V\},\qquad n\geq 2.
\]
We say a vertex algebra $V$ is $C_{n}$-\emph{cofinite} if $V/C_{n}(V)$ is finite-dimensional.
Similarly we define a vertex algebra $V$ is $E_{n}$-\emph{cofinite} if $V/E^{T}_{n}(V)$ is finite-dimensional.
\end{Definition}

Recall the following result of \cite{zhu}.

\begin{Proposition}\label{prop:zhu-poisson-algebra}
Let $V$ be a vertex algebra. Then $V/C_{2}(V)$ is {a} Poisson algebra with
\begin{align*}
&(u+C_{2}(V))\cdot(v+C_{2}(V))= u_{-1}v+C_{2}(V),\\
&[u+C_{2}(V),v+C_{2}(V)]= u_{0}v+C_{2}(V)
\end{align*}
for $u,v\in V$, and with $\mathbf{1}+C_{2}(V)$ as the identity element.
\end{Proposition}

From \cite{li-abel} and combining Definitions~\ref{dftf} and~\ref{dfCn}, we have
\begin{gather*}
E^{T}_{1}(V)=E^{T}_{2}(V)=\dots=E^{T}_{T}(V)=C_{2}(V),\\
E^{T}_{T+1}(V)=E^{T}_{T+2}(V)=\dots=E^{T}_{2T}(V)=C_{3}(V).
\end{gather*}

The following is a relation between Zhu's Poisson algebra and the degree-zero subspace $E^{T}_{0}(V)/E^{T}_{T}(V)$ of vertex Poisson algebra $\mathrm{gr}^{T}_{\mathcal{E}}(V)$.

\begin{Proposition}\label{zhu algebra}
Let $V$ be a vertex algebra and let $g$ be an automorphism of $V$ with a period $T$. Then $E^{T}_{0}(V)/E^{T}_{T}(V)$ is a Poisson algebra which coincides with Zhu's Poisson algebra $V/C_{2}(V)$, where
\begin{align*}
&\bigl(u+E^{T}_{T}(V)\bigr)\cdot\bigl(v+E^{T}_{T}(V)\bigr)= u_{-1}v+E^{T}_{T}(V),\\
&\bigl[u+E^{T}_{T}(V),v+E^{T}_{T}(V)\bigr]= u_{0}v+E^{T}_{T}(V)
\end{align*}
for $u,v \in V$.
\end{Proposition}

\section{Li filtration of twisted modules}\label{section4}

In this section, we give a construction of twisted modules for the vertex Poisson algebra $\mathrm{gr}^{T}_{\mathcal{E}}(V)$ from twisted modules {for the} vertex algebra~$V$, using a~decreasing filtration.

First, we construct the decreasing sequence $\mathcal{E}_{W}^{T}=\bigl\{E_{n}^{T}(W)\bigr\}_{n\in \mathbb{Z}}$ associated with a twisted module.

\begin{Definition}\label{dftmf} Let $V$ be a vertex algebra and let $W$ be a $g$-twisted $V$-module.
Define a~decreasing sequence $\bigl\{E_{n}^{T}(W)\bigr\}_{n\in\mathbb{Z}}$ of subspaces of $W$,
where for $n\in\mathbb{Z}$, $E_{n}^{T}(W)$ is linearly spanned by the vectors
\begin{gather}
u^{(1)}_{-1-k_{1}+\frac{r_{1}}{T}}u^{(2)}_{-1-k_{2}+\frac{r_{2}}{T}}\cdots u^{(s)}_{-1-k_{s}+\frac{r_{s}}{T}}w\label{entw-elements}
\end{gather}
for $s\geq1$, $u^{(i)}\in V^{r_{i}}$, $0\leq r_{i}\leq T-1$, $1\leq i\leq s$, $w\in W$, $k_{1}, k_{2}, \dots,k_{s}\geq 0$ with $k_{1}+k_{2}+\cdots+k_{s}-\smash{\frac{r_{1}+r_{2}+\cdots +r_{s}}{T}}\geq \frac{n}{T}$.\end{Definition}

The following are some immediate consequences.

\begin{Lemma}\label{property-entw}{Let $V$ be a vertex algebra, let $W$ be a $g$-twisted $V$-module, and let $u\in V^{r}$ with $0\leq r\leq T-1$. Then} we have
\begin{gather}
E_{n}^{T}(W)\supset E_{n+1}^{T}(W)\qquad \mbox{for $n\in\mathbb{Z}$},\nonumber
\\
E_{n}^{T}(W)=W\qquad \mbox{for $n\leq 0$},\nonumber
\\
u_{-1-k+\frac{r}{T}}E_{n}^{T}(W)\subset E_{n+kT-r}^{T}(W)\qquad \mbox{for $k\geq 0$}, \ n\in\mathbb{Z}.\label{entw-relation3}
\end{gather}\end{Lemma}

We now give a stronger spanning property for $E_{n}^{T}(W)$.

\begin{Lemma}\label{span-entw} Let $V$ be a vertex algebra and let $W$ be a $g$-twisted $V$-module. For any $n\geq1{,}$
we have
\begin{gather}
E_{n}^{T}(W)=\linspan\bigl\{u_{-1-i+\frac{r}{T}}w\mid u\in V^{r},\,i\geq1,\,w\in E_{n-iT+r}^{T}(W)\bigr\}.\label{entw-entwnitr}
\end{gather}
Furthermore, for $n\geq1$, $E_{n}^{T}(W)$ is linearly spanned by the vectors
\begin{equation}
u^{(1)}_{-1-k_{1}+\frac{r_{1}}{T}}u^{(2)}_{-1-k_{2}+\frac{r_{2}}{T}}\cdots u^{(s)}_{-1-k_{s}+\frac{r_{s}}{T}}w\label{entw-elements1}
\end{equation}
for $s\geq1$, $u^{(i)}\in V^{r_{i}}$, $0\leq r_{i}\leq T-1$, $1\leq i\leq s$, $w\in W$, $k_{1}, k_{2}, \dots, k_{s}\geq 1$ with $k_{1}+k_{2}+\cdots+k_{s}-\frac{r_{1}+r_{2}+\cdots +r_{s}}{T}\geq \frac{n}{T}$.\end{Lemma}

\begin{proof} Notice that (\ref{entw-elements1}) follows from (\ref{entw-entwnitr}) and induction. Denote by $\bigl(E_{n}^{T}\bigr)'(W)$ the subspace
on the right-hand side of (\ref{entw-entwnitr}). From (\ref{entw-relation3}), we have $\bigl(E_{n}^{T}\bigr)'(W)\subset E_{n}^{T}(W)$. In order to prove~(\ref{entw-entwnitr}), we need to prove that each spanning vector of~$E_{n}^{T}(W)$ in~(\ref{entw-elements}) lies in $\bigl(E_{n}^{T}\bigr)'(W)$. Now we use
induction on~$s$. If $s=1$, we consider $u_{-1-k+\frac{r}{T}}w\in E_{n}^{T}(W)$, $u\in V^{r}$, $0\leq r\leq T-1$, $k\geq1$. From Definition~\ref{dftmf}, we have $k-\frac{r}{T}\geq\frac{n}{T}$ such that $n-kT+r\leq0$ and $w\in W=E_{n-kT+r}^{T}(W)$,
so that we have $u_{-1-k+\frac{r}{T}}w\in \bigl(E_{n}^{T}\bigr)'(W)$. Assume $s\geq2$. If $k_{1}\geq1$, consider
\[
u^{(1)}_{-1-k_{1}+\frac{r_{1}}{T}}u^{(2)}_{-1-k_{2}+\frac{r_{2}}{T}}\cdots u^{(s)}_{-1-k_{s}+\frac{r_{s}}{T}}w\in E_{n}^{T}(W).
\]
From Definition~\ref{dftmf}, we have
\[
k_{2}+k_{3}+\cdots+k_{s}-\frac{r_{2}+r_{3}+\cdots+r_{s}}{T}\geq\frac{n-k_{1}T+r_{1}}{T}.
\]
So we get
\[
u^{(2)}_{-1-k_{2}+\frac{r_{2}}{T}}u^{(3)}_{-1-k_{3}+\frac{r_{3}}{T}}\cdots u^{(s)}_{-1-k_{s}+\frac{r_{s}}{T}}w\in E_{n-k_{1}T+r_{1}}^{T}(W).
\]
Thus we have
\[
u^{(1)}_{-1-k_{1}+\frac{r_{1}}{T}}u^{(2)}_{-1-k_{2}+\frac{r_{2}}{T}}\cdots u^{(s)}_{-1-k_{s}+\frac{r_{s}}{T}}w\in \bigl(E_{n}^{T}\bigr)'(W).
\]
If $k_{1}=0$, from Definition~\ref{dftmf} we have
\[
u^{(2)}_{-1-k_{2}+\frac{r_{2}}{T}}u^{(3)}_{-1-k_{3}+\frac{r_{3}}{T}}\cdots u^{(s)}_{-1-k_{s}+\frac{r_{s}}{T}}w\in E_{n+r_{1}}^{T}(W).
\]
By the inductive hypothesis, we have
\[
u^{(2)}_{-1-k_{2}+\frac{r_{2}}{T}}u^{(3)}_{-1-k_{3}+\frac{r_{3}}{T}}\cdots u^{(s)}_{-1-k_{s}+\frac{r_{s}}{T}}w\in \bigl(E_{n+r_{1}}^{T}\bigr)'(W).
\]
 Furthermore, for any $b\in V^{p}$, $0\leq p\leq T-1$, $k\geq 1$, $w'\in E_{n+r_{1}-kT+p}^{T}(W)$, we have
\[
u^{(1)}_{-1+\frac{r_{1}}{T}}b_{-1-k+\frac{p}{T}}w'=b_{-1-k+\frac{p}{T}}u^{(1)}_{-1+\frac{r_{1}}{T}}w'+\sum_{i\in \mathbb{N}}\binom{-1+\frac{r_{1}}{T}}{i}
(u^{(1)}_{i}b)_{-2-k-i+\frac{r_{1}+p}{T}}w'.
\]
From Definition~\ref{dftmf} we have $u^{(1)}_{-1+\frac{r_{1}}{T}}w'\in E_{n-kT+p}^{T}(W)$, so that $b_{-1-k+\frac{p}{T}}u^{(1)}_{-1+\frac{r_{1}}{T}}w'\in \bigl(E_{n}^{T}\bigr)'(W)$.
On the other hand, if $0\leq r_{1}+p\leq T-1$, we have \[
\bigl(u^{(1)}_{i}b\bigr)_{-2-k-i+\frac{r_{1}+p}{T}}w'=(u^{(1)}_{i}b)_{-1-(k+i+1)+\frac{r_{1}+p}{T}}w',
\] with
$w'\in E_{n-kT+r_{1}+p}^{T}(W)\subset E_{n-(1+k+i)T+r_{1}+p}^{T}(W)$.
So that
\[
\bigl(u^{(1)}_{i}b\bigr)_{-2-k-i+\frac{r_{1}+p}{T}}w'\in \bigl(E_{n}^{T}\bigr)'(W).
\]
If $T\leq r_{1}+p \leq 2T-2$, we have \[
\bigl(u^{(1)}_{i}b\bigr)_{-2-k-i+\frac{r_{1}+p}{T}}w'=\bigl(u^{(1)}_{i}b\bigr)_{-1-(k+i)+\frac{r_{1}+p-T}{T}}w'
\]
with $w'\in E_{n-kT+r_{1}+p}^{T}(W)\subset E_{n-(k+i)T+r_{1}+p-T}^{T}(W)$, so that
\[
\bigl(u^{(1)}_{i}b\bigr)_{-2-k-i+\frac{r_{1}+p}{T}}w'\in \bigl(E_{n}^{T}\bigr)'(W).
\]
Therefore, $u^{(1)}_{-1+\frac{r_{1}}{T}}b_{-1-k+\frac{p}{T}}w'\in \bigl(E_{n}^{T}\bigr)'(W)$, this proves that
\[
u^{(1)}_{-1+\frac{r_{1}}{T}}u^{(2)}_{-1-k_{2}+\frac{r_{2}}{T}}\cdots u^{(s)}_{-1-k_{s}+\frac{r_{s}}{T}}w\in \bigl(E_{n}^{T}\bigr)'(W),
\]
completing the proof.
\end{proof}

\begin{Lemma}\label{am-entw} Let $V$ be a vertex algebra and let $W$ be a $g$-twisted $V$-module. For any $a\in V^{p}$, $0\leq p \leq T-1$, $m,n\in\mathbb{Z}$, we have
\begin{gather}
a_{m+\frac{p}{T}}E_{n}^{T}(W)\subset E_{n-(m+1)T-p}^{T}(W)\qquad \mbox{for $m\in \mathbb{Z}$.}\label{a-entw1}
\end{gather}
Furthermore,
\begin{gather}
a_{m+\frac{p}{T}}E_{n}^{T}(W)\subset E_{n-(m+1)T-p+T}^{T}(W)\qquad \mbox{for $m\geq 0$}.\label{a-entw2}
\end{gather}
\end{Lemma}

\begin{proof}If $m\leq-1$, by (\ref{entw-relation3}) we have
\[
a_{m+\frac{p}{T}}E_{n}^{T}(W)=a_{-1-(-1-m)+\frac{p}{T}}E_{n}^{T}(W)\subset E_{n-(m+1)T-p}^{T}(W).
\]
 Assume $m\geq0$, since $E_{n-(m+1)T-p+T}^{T}(W)\subset E_{n-(m+1)T-p}^{T}(W)$, it suffices to prove (\ref{a-entw2}). We now prove the assertion by induction on $n$. If $n\leq0$, we have $n-(m+1)T-p+T\leq0$, we have $a_{m+\frac{p}{T}}E_{n}^{T}(W)\subset W=E_{n-(m+1)T-p+T}^{T}(W).$ If $m\geq 0$, $n\geq 1$, from (\ref{entw-entwnitr}) we have $E_{n}^{T}(W)$ is spanned by the vectors $u_{-1-k+\frac{q}{T}}w$ with $u\in V^{q}$, $0\leq q \leq T-1$, $k\geq 1$, $w\in E_{n-kT+q}^{T}(W)$. In view of~(\ref{twisted-com}), we have
\[
a_{m+\frac{p}{T}}u_{-1-k+\frac{q}{T}}w=u_{-1-k+\frac{q}{T}}a_{m+\frac{p}{T}}w+\sum_{i\in \mathbb{N}}\binom{m+\frac{p}{T}}{i}(a_{i}u)_{m-k-i-1+\frac{p+q}{T}}w.
\]
Since $w\in E_{n-kT+q}^{T}(W)$, $n-kT+q<n$, from the inductive hypothesis and (\ref{entw-relation3}) we have
\begin{gather*}
a_{m+\frac{p}{T}}w\in E_{n-kT+q-(m+1)T-p+T}^{T}(W),\\
u_{-1-k+\frac{q}{T}}a_{m+\frac{p}{T}}w\in u_{-1-k+\frac{q}{T}}E_{n-kT+q-(m+1)T-p+T}^{T}(W)\subset E_{n-(m+1)T-p+T}^{T}(W).
\end{gather*}
If $0\leq p+q\leq T-1$, since $i\geq0$, from the inductive hypothesis we have
\[
(a_{i}u)_{m-k-i-1+\frac{p+q}{T}}w\in E_{n-kT+q-(m-k-i-1+1)T-p-q}^{T}(W)\subset E_{n-(m+1)T-p+T}^{T}(W).
\]
If $T\leq p+q\leq 2T-2$, since $i\geq 0$, similarly we have
\begin{align*}
\begin{split}
& (a_{i}u)_{m-k-i-1+\frac{p+q}{T}}w\\
&\qquad {}=(a_{i}u)_{m-k-i+\frac{p+q-T}{T}}w \in E_{n-kT+q-(m-k-i+1)T-p-q+T}^{T}(W)
 \subset E_{n-(m+1)T-p+T}^{T}(W).
\end{split}
\end{align*}
Therefore, $a_{m+\frac{p}{T}}u_{-1-k+\frac{q}{T}}w\in E_{n-(m+1)T-p+T}^{T}(W)$. This proves
\[
a_{m+\frac{p}{T}}E_{n}^{T}(W)\subset E_{n-(m+1)T-p+T}^{T}(W),
\]
 completing the induction and the whole proof.
\end{proof}

Now we have the following general case.

\begin{Proposition}\label{unw-entw} Suppose that $V$ is a vertex algebra, let $W$ be a $g$-twisted $V$-module and let $u\in E_{rT}^{T}(V)\cap V^{p}$, $w\in E_{s}^{T}(W)$ with $0\leq p \leq T-1$, $r,s\in \mathbb{Z}$. Then,
\begin{gather}
 u_{n+\frac{p}{T}}w\in E_{rT+s-(n+1)T-p}^{T}(W)\qquad \mbox{for $n\in\mathbb{Z}$}.\label{general-tmrs3}
\end{gather}
\end{Proposition}

\begin{proof}We are going to use induction on $r$. By~(\ref{a-entw1}), we have $u_{n+\frac{p}{T}}w\in E_{s-(n+1)T-p}^{T}(W)$.
If $r\leq0$, we have $rT+s-(n+1)T-p\leq s-(n+1)T-p$, so that $u_{n+\frac{p}{T}}w\in E_{s-(n+1)T-p}^{T}(W)\subset E_{rT+s-(n+1)T-p}^{T}(W)$.

Assume $r\geq 1$ and \smash{$u\in E_{(r+1)T}^{T}(V)$}, $u\in V^{p}$, $0\leq p \leq T-1$. In view of~(\ref{spann-i}), it suffices to consider $u=a_{-2-k}b$ with $a\in V^{q}$, $b\in V^{t}$, \smash{$b\in E_{(r-k)T}^{T}(V)$}, $0\leq k\leq r$, $0\leq q, t\leq T-1$. And from~(\ref{twisted-weakcom}), we know that for each $a\in V^{q}$ and $b\in V^{t}$, there exists a nonnegative integer $N$ such that $(x_{1}-x_{2})^{N}[Y_{W}(a,x_{1}),Y_{W}(b,x_{2})]w=0$.

If $0\leq q+t \leq T-1$, then $q+t=p$. By (\ref{twisted-ite-1}), we have
\begin{align*}
(a_{-2-k}b)_{n+\frac{q+t}{T}}w &= \sum_{i\in \mathbb{N}}\sum_{-2-k\leq j<N}((-1)^{i}\binom{-\frac{q}{T}}{j+k+2}\binom{j}{i}a_{j-i+\frac{q}{T}}b_{n-2-k+i-j+\frac{t}{T}}w\\
&\quad -(-1)^{i+j}\binom{-\frac{q}{T}}{j+k+2}\binom{j}{i}b_{n-2-k-i+\frac{t}{T}}a_{i+\frac{q}{T}}w).
\end{align*}
For $n\in\mathbb{Z}$, $b\in E_{(r-k)T}^{T}(V)$ with $r-k\leq r$, using the inductive hypothesis and (\ref{a-entw1}) we have
\begin{align*}
&a_{j-i+\frac{q}{T}}b_{n-2-k+i-j+\frac{t}{T}}w\in a_{j-i+\frac{q}{T}}E^{T}_{(r-k)T+s-(n-2-k+i-j+1)T-t}(W)\\
&\hphantom{a_{j-i+\frac{q}{T}}b_{n-2-k+i-j+\frac{t}{T}}w}{}
\subset E^{T}_{(r-k)T+s-(n-2-k+i-j+1)T-t-(j-i+1)T-q}(W)\\
&\hphantom{a_{j-i+\frac{q}{T}}b_{n-2-k+i-j+\frac{t}{T}}w}{}
=E^{T}_{(r+1)T+s-(n+1)T-t-q}(W)
=E^{T}_{(r+1)T+s-(n+1)T-p}(W),\\
&b_{n-2-k-i+\frac{t}{T}}a_{i+\frac{q}{T}}w\in b_{n-2-k-i+\frac{t}{T}}E^{T}_{s-(i+1)T-q}(W)\\
&\hphantom{b_{n-2-k-i+\frac{t}{T}}a_{i+\frac{q}{T}}w}{}
\subset E^{T}_{(r-k)T+s-(i+1)T-q-(n-2-k-i+1)T-t}(W)\\
&\hphantom{b_{n-2-k-i+\frac{t}{T}}a_{i+\frac{q}{T}}w}{}
=E^{T}_{(r+1)T+s-(n+1)T-q-t}(W)
=E^{T}_{(r+1)T+s-(n+1)T-p}(W),
\end{align*}
from which we have $u_{n+\frac{p}{T}}w\in E_{(r+1)T+s-(n+1)T-p}^{T}(W)$.

If $T\leq q+t \leq 2T-2$, then $q+t-T=p$. By (\ref{twisted-ite-1}), we have
\begin{align*}
(a_{-2-k}b)_{n+\frac{q+t-T}{T}}w &= \sum_{i\in \mathbb{N}}\sum_{-2-k\leq j<N}((-1)^{i}\binom{-\frac{q}{T}}{j+k+2}\binom{j}{i}a_{j-i+\frac{q}{T}}b_{n-3-k+i-j+\frac{t}{T}}w\\
&\quad -(-1)^{i+j}\binom{-\frac{q}{T}}{j+k+2}\binom{j}{i}b_{n-3-k-i+\frac{t}{T}}a_{i+\frac{q}{T}}w).
\end{align*}
Similarly, using inductive hypothesis and (\ref{a-entw1}), we have
\begin{align*}
&a_{j-i+\frac{q}{T}}b_{n-3-k+i-j+\frac{t}{T}}w\in a_{j-i+\frac{q}{T}}E^{T}_{(r-k)T+s-(n-3-k+i-j+1)T-t}(W)\\
&\hphantom{a_{j-i+\frac{q}{T}}b_{n-3-k+i-j+\frac{t}{T}}w}{}
\subset E^{T}_{(r-k)T+s-(n-3-k+i-j+1)T-t-(j-i+1)T-q}(W)\\
&\hphantom{a_{j-i+\frac{q}{T}}b_{n-3-k+i-j+\frac{t}{T}}w}{}=E^{T}_{(r+1)T+s-(n+1)T+T-t-q}(W)
 =E^{T}_{(r+1)T+s-(n+1)T-p}(W),\\
&b_{n-3-k-i+\frac{t}{T}}a_{i+\frac{q}{T}}w\in b_{n-3-k-i+\frac{t}{T}}E^{T}_{s-(i+1)T-q}(W)\\
&\hphantom{b_{n-3-k-i+\frac{t}{T}}a_{i+\frac{q}{T}}w}{}
\subset E^{T}_{(r-k)T+s-(i+1)T-q-(n-3-k-i+1)T-t}(W)\\
&\hphantom{b_{n-3-k-i+\frac{t}{T}}a_{i+\frac{q}{T}}w}{} =E^{T}_{(r+1)T+s-(n+1)T+T-q-t}(W)
 =E^{T}_{(r+1)T+s-(n+1)T-p}(W),
\end{align*}
from which we have $u_{n+\frac{p}{T}}w\in E_{(r+1)T+s-(n+1)T-p}^{T}(W)$. This concludes the proof.
\end{proof}

The case $T=1$ is Li's result \cite{li-abel}.
The following is the main result of this section.

\begin{Proposition}\label{pgretmodule} {Let $V$ be a vertex algebra, let $g$ have finite period $T$, and let} $W$ be a $g$-twisted $V$-module. For $T\geq2$, let $E^{T}_{n}(W)$ be the decreasing sequence defined in Definition~{\rm \ref{dftmf}} for $W$,
then the associated graded vector space $\mathrm{gr}^{T}_{\mathcal{E}}(W)=\coprod_{n\geq0}E_{n}^{T}(W)/E_{n+1}^{T}(W)$ is a $g$-twisted module for the vertex Poisson algebra $\mathrm{gr}^{T}_{\mathcal{E}}(V)$ with
\begin{gather}
\bigl(u+E_{(r+1)T}^{T}(V)\bigr)\cdot\bigl(w+E_{s+1}^{T}(W)\bigr)
=u_{-1+\frac{p}{T}}w+E_{rT+s-p+1}^{T}(W),\label{pgretmodule1}\\
Y_{-}^{W}(u+E_{(r+1)T}^{T}(V), x)\bigl(w+E_{s+1}^{T}(W)\bigr)\nonumber\\
\qquad = \begin{cases} \displaystyle \sum_{n\in \mathbb{N}}\bigl(u_{n}w+E_{rT+s-(n+1)T+1}^{T}(W)\bigr)x^{-n-1} & \mbox{for $p=0$},\\
\displaystyle \sum_{n\in \mathbb{Z}}\bigl(u_{n+\frac{p}{T}}w+E_{rT+s-(n+1)T-p+1}^{T}(W)\bigr)x^{-n-\frac{p}{T}-1}& \mbox{for $1\leq p\leq T-1$},\end{cases}\label{pgretmodule2}
\end{gather}
with $u\in E_{rT}^{T}(V)\cap V^{p}$, $w\in E_{s}^{T}(W)$, $r,s\in\mathbb{N}$, $0\leq p\leq T-1$.
\end{Proposition}

\begin{proof}
From Proposition~\ref{unw-entw}, the actions given by (\ref{pgretmodule1}) and (\ref{pgretmodule2}) are well defined. The vacuum axiom gives $\mathbf{1}_{-1}w=w$, hence $\bigl(\mathbf{1}+E_T^T(V)\bigr)\cdot\bigl(w+E_{s+1}^T(W)\bigr)=w+E_{s+1}^T(W)$. Clearly,
it is straightforward to check~(\ref{dftmvla-1}) and~(\ref{dftmvla-2}). Now we prove the compatibility~(\ref{comp}).
For $u\in E^{T}_{rT}{(V)}\cap V^{p}$, $v\in E^{T}_{sT}{(V)}\cap V^{q}$, $w\in E^{T}_{k}(W)$, $r,s,{k}\in\mathbb{N}$.

For $p=0$, $0\leq q \leq T-1$, by (\ref{pgretmodule2}) we have $n\in\mathbb{N}$. From (\ref{twisted-com}), we have
\[
u_{n}v_{-1+\frac{q}{T}}w=v_{-1+\frac{q}{T}}u_{n}w+\sum_{i=0}^{n}\binom{n}{i}(u_{i}v)_{n-1-i+\frac{q}{T}}w.
\]
Since $i\geq0$, by (\ref{general-tmrs3}) and (\ref{general-rs0}) we have
\begin{gather*}
(u_{i}v)_{n-1-i+\frac{q}{T}}w\in E^{T}_{(r+s-n)T+k-q}(W)\subset E^{T}_{(r+s-n-1)T+k-q+1}(W),\\
u_{n}v_{-1+\frac{q}{T}}w\in E^{T}_{(r+s-n-1)T+k-q}(W),\\
v_{-1+\frac{q}{T}}u_{n}w\in E^{T}_{(r+s-n-1)T+k-q}(W).
\end{gather*}
Then
\[
u_{n}(v_{-1+\frac{q}{T}}w)=v_{-1+\frac{q}{T}}(u_{n}w)+(u_{n}v)_{-1+\frac{q}{T}}w+\sum_{i=0}^{n-1}\binom{n}{i}(u_{i}v)_{n-1-i+\frac{q}{T}}w.
\]
Thus
\begin{align*}
& u_{n}(v_{-1+\frac{q}{T}}w)+E^{T}_{(r+s-n-1)T+k-q+1}(W)\\
&\qquad = v_{-1+\frac{q}{T}}(u_{n}w)+(u_{n}v)_{-1+\frac{q}{T}}w+E^{T}_{(r+s-n-1)T+k-q+1}(W).
\end{align*}

For $1\leq p\leq T-1$, $0\leq q \leq T-1$, then by (\ref{pgretmodule2}) we have $n\in\mathbb{Z}$.
From (\ref{twisted-com}), we have
\[
u_{n+\frac{p}{T}}v_{-1+\frac{q}{T}}w=v_{-1+\frac{q}{T}}u_{n+\frac{p}{T}}w+\sum_{i\in \mathbb{N}}\binom{n+\frac{p}{T}}{i}(u_{i}v)_{n-1-i+\frac{p+q}{T}}w.
\]
For $n\geq0$, since $i\geq0$, by (\ref{general-tmrs3}) and (\ref{general-rs0}) we have
\begin{gather*}
(u_{i}v)_{n-1-i+\frac{p+q}{T}}w\in E^{T}_{(r+s-n)T+k-p-q}(W)\subset E^{T}_{(r+s-n-1)T+k-p-q+1}(W),\\
u_{n+\frac{p}{T}}v_{-1+\frac{q}{T}}w\in E^{T}_{(r+s-n-1)T+k-p-q}(W),\\
v_{-1+\frac{q}{T}}u_{n+\frac{p}{T}}w\in E^{T}_{(r+s-n-1)T+k-p-q}(W).
\end{gather*}
If $1\leq p+q\leq T-1$, we have
\[
(u_{n}v)_{-1+\frac{p+q}{T}}w\in E^{T}_{(r+s-n)T+k-p-q}(W)\subset E^{T}_{(r+s-n-1)T+k-p-q+1}(W).
\]
Then
\begin{align*}
& u_{n+\frac{p}{T}}(v_{-1+\frac{q}{T}}w)+E^{T}_{(r+s-n)T+k-p-q+1}(W)\\
&\qquad = (u_{n}v)_{-1+\frac{p+q}{T}}w+v_{-1+\frac{q}{T}}(u_{n+\frac{p}{T}}w)+E^{T}_{(r+s-n)T+k-p-q+1}(W).
\end{align*}
If $T\leq p+q\leq 2T-2$, we also have
\[
(u_{n}v)_{-1+\frac{p+q-T}{T}}w\in E^{T}_{(r+s-n)T+k-p-q}(W)\subset E^{T}_{(r+s-n-1)T+k-p-q+1}(W).
\]
Then
\begin{align*}
& u_{n+\frac{p}{T}}(v_{-1+\frac{q}{T}}w)+E^{T}_{(r+s-n)T+k-p-q+1}(W)\\
&\qquad = (u_{n}v)_{-1+\frac{p+q-T}{T}}w+v_{-1+\frac{q}{T}}(u_{n+\frac{p}{T}}w)+E^{T}_{(r+s-n)T+k-p-q+1}(W).
\end{align*}
For $n\leq-1$, by (\ref{twisted-com}) we have
\[
u_{n+\frac{p}{T}}v_{-1+\frac{q}{T}}w=v_{-1+\frac{q}{T}}u_{n+\frac{{p}}{T}}w+\sum_{i\in \mathbb{N}}\binom{n+\frac{p}{T}}{i}(u_{i}v)_{n-1-i+\frac{p+q}{T}}w.
\]
Similarly, we have
\begin{gather*}
(u_{i}v)_{n-1-i+\frac{p+q}{T}}w\in E^{T}_{(r+s-n)T+k-p-q}(W)\subset E^{T}_{(r+s-n-1)T+k-p-q+1}(W),\\
u_{n+\frac{p}{T}}v_{-1+\frac{q}{T}}w\in E^{T}_{(r+s-n-1)T+k-p-q}(W),\\
v_{-1+\frac{q}{T}}u_{n+\frac{p}{T}}w\in E^{T}_{(r+s-n-1)T+k-p-q}(W).
\end{gather*}
Then
\begin{align*}
& u_{n+\frac{p}{T}}(v_{-1+\frac{q}{T}}w)+E^{T}_{(r+s-n)T+k-p-q+1}(W)\\
&\qquad = {(u_{n}v)_{-1+\frac{p+q}{T}}w+}v_{-1+\frac{q}{T}}(u_{n+\frac{p}{T}}w)+E^{T}_{(r+s-n)T+k-p-q+1}(W).
\end{align*}
To sum up, it proves $Y^{W}_{-}(u,x)(vw)=(Y_{-}(u,x)v)w+vY^{W}_{-}(u,x)w$.
Therefore $\mathrm{gr}^{T}_{\mathcal{E}}(W)$ is a~$g$\nobreakdash-twisted module for the vertex Poisson algebra $\mathrm{gr}^{T}_{\mathcal{E}}(V)$.
\end{proof}

In Section~\ref{section6}, we will {exclude} the possibility that the associated sequence $\mathcal{E}^{T}_{W}$ is trivial in the sense that $W=E^{T}_{n}(W)$ for all $n\geq0$. Here we have the following lemma.

\begin{Lemma}\label{wt-entw}
Let $V=\coprod_{n\geq 0}V_{(n)}$ be an $\mathbb{N}$-graded vertex algebra, let \smash{$W=\bigoplus_{n\in \frac{1}{T}\mathbb{N}}W_{(n)}$} be a~$\frac{1}{T}\mathbb{N}$-graded $g$-twisted $V$-module and $\mathcal{E}^{T}_{W}=\bigl\{E_{n}^{T}(W)\bigr\}_{n \in \mathbb{Z}}$ be the decreasing sequence defined in Definition~{\rm \ref{dftmf}}. Then
\begin{gather}
E^{T}_{n}(W)\subset \coprod_{m\geq\frac{n}{T}}W_{(m)}\qquad \mbox{for $n\geq 0$}.\label{wt-entw-1}
\end{gather}
Furthermore, the associated decreasing sequence $\mathcal{E}^{T}_{W}=\bigl\{E_{n}^{T}(W)\bigr\}_{n \in \mathbb{Z}}$ for~$W$ is a filtration, i.e.,
\begin{gather}
\bigcap_{n\geq 0}E^{T}_{n}(W)=0.\label{wt-entw-2}
\end{gather}
\end{Lemma}

\begin{proof}
For $n=0$, from Definition~\ref{dftmf} we have $E^{T}_{0}(W)=\coprod_{m\geq0}W_{(m)}$. For $ n\geq 1$, $E^{T}_{n}(W)$ is linearly spanned by the vectors
\[
u^{(1)}_{-1-k_{1}+\frac{r_{1}}{T}}u^{(2)}_{-1-k_{2}+\frac{r_{2}}{T}}\cdots u^{(s)}_{-1-k_{s}+\frac{r_{s}}{T}}w
\]
 with $s\geq 1$, $u^{(i)}\in V^{r_{i}}$, $w\in W_{(m_{1})}$, $k_{1}, k_{2}, \dots, k_{s}\geq 1$, $k_{1}+k_{2}+\cdots + k_{s}-\frac{r_{1}+r_{2}+\cdots+r_{s}}{T}\geq\frac{n}{T}$, $m_{1}\geq 0$, $0\leq r_{i}\leq T-1$, $1\leq i\leq s$.
If the vectors $u^{(1)}, u^{(2)}, \dots,u^{(s)},w$ are homogeneous, we have
\begin{align*}
& \wt (u^{(1)}_{-1-k_{1}+\frac{r_{1}}{T}}u^{(2)}_{-1-k_{2}+\frac{r_{2}}{T}}\cdots u^{(s)}_{-1-k_{s}+\frac{r_{s}}{T}}w)\\
&\quad = \wt u^{(1)}+\wt u^{(2)}+\cdots + \wt u^{(s)}+k_{1}+k_{2}+\cdots+k_{s}-\frac{r_{1}+r_{2}+\cdots+r_{s}}{T}+m_{1}
 \geq \frac{n}{T}.
\end{align*}
This proves (\ref{wt-entw-1}).

Clearly, each subspace $E^{T}_{n}(W)$ of $W$ is graded, we immediately have (\ref{wt-entw-2}).
\end{proof}

\section[The relation between the sequences E(W,T) and C(W,T)]{The relation between the sequences $\boldsymbol{\mathcal{E}^{T}_{W}}$ and $\boldsymbol{\mathcal{C}^{T}_{W}}$}\label{section5}

In this section, we introduce the decreasing sequence \smash{$\mathcal{C}^{T}_{W}=\bigl\{C_{n}^{T}(W)\bigr\}_{n\in \mathbb{Z}_{\geq 2}}$} of twisted modules.
We define the subspace $C^{T}_{n}(W)$ of a twisted module $W$ as follows.

\begin{Definition}\label{dftmcn}
Let $V$ be a vertex algebra and let $W$ be a $g$-twisted $V$-module. For $n\geq 2$, we define $C^{T}_{n}(W)$ to be the subspace of~$W$, with
\begin{gather}
C^{T}_{n}(W)=\linspan\bigl\{u_{-n+\frac{p}{T}}w \mid u\in V^{p},\, w\in W,\, p=0, 1, \dots, T-1\bigr\}.\label{cntw-def}
\end{gather}
\end{Definition}

A twisted $V$-module $W$ is said to be \emph{$C_{n}$-cofinite} if $W/C^{T}_{n}(W)$ is finite-dimensional.
A twisted $V$-module $W$ is said to be \emph{$E_{n}$-cofinite} if $W/E^{T}_{n}(W)$ is finite-dimensional.

Here are some consequences for $C^{T}_{n}(W)$.
\begin{Lemma}\label{property-cntw}Let $V$ be a vertex algebra and let $W$ be a $g$-twisted $V$-module. For $n\geq 2$, we have
\begin{gather}
C^{T}_{m}(W)\subset C^{T}_{n}(W)\ \ \mbox{for $m\geq n$}, \label{tmcn-relation-1}\\
u_{-k+\frac{p}{T}}C^{T}_{n}(W)\subset C^{T}_{n}(W), \label{tmcn-relation-2}\\
u_{-n+\frac{p}{T}}v_{-k+\frac{q}{T}}w\equiv v_{-k+\frac{q}{T}}u_{-n+\frac{p}{T}}w\ \bigl(\mod\ C^{T}_{n+k-1}(W)\bigr)\label{tmcn-relation-3}
\end{gather}
for $u\in V^{p}$, $v\in V^{q}$, $0\leq p, q \leq T-1$, $w\in W$, $k\geq 1$.
\end{Lemma}

\begin{proof}
For $u\in V^{p}$, $0\leq p \leq T-1$, $w\in W$, $r\geq 2$, we have
\[
u_{-r-1+\frac{p}{T}}w=\frac{1}{r-\frac{p}{T}}(\D u)_{-r+\frac{p}{T}}w.
\] From this, we immediately have $C^{T}_{r+1}(W)\subset C^{T}_{r}(W)$ for $r\geq2$, this implies~(\ref{tmcn-relation-1}).

For $u\in V^{p}$, $v\in V^{q}$, $0\leq p, q \leq T-1$, $w\in W$, $k\geq 1$, $n\geq 2$, by~(\ref{twisted-com}) we have
\[
u_{-k+\frac{p}{T}}v_{-n+\frac{q}{T}}w=v_{-n+\frac{q}{T}}u_{-k+\frac{p}{T}}w+\sum_{i\in \mathbb{N}}\binom{-k+\frac{p}{T}}{i}(u_{i}v)_{-n-k-i+\frac{p+q}{T}}w.
\]
If $0\leq p+q\leq T-1$, then
\[
(u_{i}v)_{-n-k-i+\frac{p+q}{T}}w\in C^{T}_{n+k+i}(W)\subset C^{T}_{n}(W).
\]
 If $T\leq p+q\leq 2T-2$, then we have
 \[
 (u_{i}v)_{-n-k-i+1+\frac{p+q-T}{T}}w\in C^{T}_{n+k+i-1}(W)\subset C^{T}_{n}(W),
 \] this proves~(\ref{tmcn-relation-2}).

We also have
\begin{gather}
u_{-n+\frac{p}{T}}v_{-k+\frac{q}{T}}w-v_{-k+\frac{q}{T}}u_{-n+\frac{p}{T}}w=\sum_{i\in \mathbb{N}}\binom{-n+\frac{p}{T}}{i}(u_{i}v)_{-n-k-i+\frac{p+q}{T}}w.\label{property-cntw-1}
\end{gather}
If $0\leq p+q \leq T-1$, we have
\begin{gather}
\sum_{i\in \mathbb{N}}\binom{-n+\frac{p}{T}}{i}(u_{i}v)_{-n-k-i+\frac{p+q}{T}}w\in C_{n+k}^{T}(W)\subset C^{T}_{n+k-1}(W),\label{property-cntw-2}
\end{gather}
and if $T\leq p+q\leq 2T-2$, we have
\begin{gather}
 \sum_{i\in \mathbb{N}}\binom{-n+\frac{p}{T}}{i}(u_{i}v)_{-n-k-i+\frac{p+q}{T}}w\nonumber\\
\qquad =\sum_{i\in \mathbb{N}}\binom{-n+\frac{p}{T}}{i}(u_{i}v)_{-n-k-i+1+\frac{p+q-T}{T}}w\in C_{n+k-1}^{T}(W).\label{property-cntw-3}
\end{gather}
Combining (\ref{property-cntw-1}), (\ref{property-cntw-2}) and (\ref{property-cntw-3}), we have (\ref{tmcn-relation-3}).
\end{proof}

We have the following technical result.

\begin{Lemma}\label{u(T+1)C3}
Let $V$ be a vertex algebra and let $W$ be a $g$-twisted $V$-module. Then
\[
u^{(1)}_{-2+\frac{r_{1}}{T}}u^{(2)}_{-2+\frac{r_{2}}{T}}\cdots u^{(s)}_{-2+\frac{r_{s}}{T}}w\in C^{T}_{3}(W)
\]
for $s\geq T+1$, $u^{(i)}\in V^{r_{i}}$, $w\in W$, $0\leq r_{i}\leq T-1$, $1\leq i\leq s$.
\end{Lemma}
\begin{proof}
For $u\in V^{p}$, $0\leq p\leq T-1$, since $u_{-2+\frac{p}{T}}C^{T}_{3}(W)\subset C^{T}_{3}(W)$, it is sufficient to prove the conclusion holds for $s=T+1$. Assume $s=T+1$, if there exist $k_{1}$, $k_{2}$ with $1\leq k_{1}< k_{2}\leq s$ such that $0\leq r_{k_{1}}+r_{k_{2}}\leq T-1$, then by~(\ref{twisted-ite-2}), there exist nonnegative integers $k$ and $l$ such that
\begin{align*}
& \bigl(u^{(k_{1})}_{-1}u^{(k_{2})}\bigr)_{-3+\frac{r_{k_{1}}+r_{k_{2}}}{T}}w\\
&\qquad = \sum_{i=0}^{k}\sum_{j\in \mathbb{N}}\binom{-l-\frac{r_{k_{1}}}{T}}{i}\binom{i-1}{j}(-1)^{j}u^{(k_{1})}_{-1+l+i-j+\frac{r_{k_{1}}}{T}}u^{(k_{2})}_{-3-l-i+j+\frac{r_{k_{2}}}{T}}w.
\end{align*}

Since $0\leq r_{k_{1}}+r_{k_{2}}\leq T-1$, we have \smash{$\bigl(u^{(k_{1})}_{-1}u^{(k_{2})}\bigr)_{-3+\frac{r_{k_{1}}+r_{k_{2}}}{T}}w\in C^{T}_{3}(W)$}.
If $i-1\geq 0$, then $j\leq i-1$. For each $i$, $j$, from~(\ref{twisted-com}) we have
\begin{align*}
& u^{(k_{1})}_{-1+l+i-j+\frac{r_{k_{1}}}{T}}u^{(k_{2})}_{-3-l-i+j+\frac{r_{k_{2}}}{T}}w\\
& = u^{(k_{2})}_{-3-l-i+j+\frac{r_{k_{2}}}{T}}u^{(k_{1})}_{-1+l+i-j+\frac{r_{k_{1}}}{T}}w
 +\sum_{n\in \mathbb{N}}\binom{-1+l+i-j+\frac{r_{k_{1}}}{T}}{n}(u^{(k_{1})}_{n}u^{(k_{2})})_{-4-n+\frac{r_{k_{1}}+r_{k_{2}}}{T}}w.
\end{align*}
Notice that $-4-n\leq -3$ and $-3-l-i+j\leq -3$, we have
\begin{gather*}
\bigl(u^{(k_{1})}_{n}u^{(k_{2})}\bigr)_{-4-n+\frac{r_{k_{1}}+r_{k_{2}}}{T}}w\in C^{T}_{3}(W),\\
u^{(k_{2})}_{-3-l-i+j+\frac{r_{k_{2}}}{T}}u^{(k_{1})}_{-1+l+i-j+\frac{r_{k_{1}}}{T}}w\in C^{T}_{3}(W),
\end{gather*}
which implies
\[
u^{(k_{1})}_{-1+l+i-j+\frac{r_{k_{1}}}{T}}u^{(k_{2})}_{-3-l-i+j+\frac{r_{k_{2}}}{T}}w\in C^{T}_{3}(W).
\]

For $i=0$, if $j\geq l+2$, we have $-1+l-j\leq -3$, so that \[
u^{(k_{1})}_{-1+l-j+\frac{r_{k_{1}}}{T}}u^{(k_{2})}_{-3-l+j+\frac{r_{k_{2}}}{T}}w\in C^{T}_{3}(W).
\]
If $j\leq l$, we have $-3-l+j\leq-3$. From~(\ref{twisted-com}), we have
\begin{align*}
& u^{(k_{1})}_{-1+l-j+\frac{r_{k_{1}}}{T}}u^{(k_{2})}_{-3-l+j+\frac{r_{k_{2}}}{T}}w\\
&\qquad = u^{(k_{2})}_{-3-l+j+\frac{r_{k_{2}}}{T}}u^{(k_{1})}_{-1+l-j+\frac{r_{k_{1}}}{T}}w+\sum_{n\in \mathbb{N}}\binom{-1+l-j+\frac{r_{k_{1}}}{T}}{n}\bigl(u^{(k_{1})}_{n}u^{(k_{2})}\bigr)_{-4-n+\frac{r_{k_{1}}+r_{k_{2}}}{T}}w.
\end{align*}
Thus we have
\[
u^{(k_{1})}_{-1+l-j+\frac{r_{k_{1}}}{T}}u^{(k_{2})}_{-3-l+j+\frac{r_{k_{2}}}{T}}w\in C^{T}_{3}(W).
\] Therefore, the only remaining case $i=0$ and $j=l+1$, the corresponding term \smash{$u^{(k_{1})}_{-2+\frac{r_{k_{1}}}{T}}u^{(k_{2})}_{-2+\frac{r_{k_{2}}}{T}}w$} must also lie in $C_{3}^{T}(W)$. For any $w\in W$, by (\ref{tmcn-relation-3}) we have
\begin{align*}
& u^{(1)}_{-2+\frac{r_{1}}{T}}u^{(2)}_{-2+\frac{r_{2}}{T}}\cdots u^{(s)}_{-2+\frac{r_{s}}{T}}w
 \equiv u^{(k_{1})}_{-2+\frac{r_{k_{1}}}{T}}u^{(k_{2})}_{-2+\frac{r_{k_{2}}}{T}}u^{(1)}_{-2+\frac{r_{1}}{T}}\cdots u^{(k_{1}-1)}_{-2+\frac{r_{k_{1}-1}}{T}}u^{(k_{1}+1)}_{-2+\frac{r_{k_{1}+1}}{T}}\cdots\\
&\hphantom{u^{(1)}_{-2+\frac{r_{1}}{T}}u^{(2)}_{-2+\frac{r_{2}}{T}}\cdots u^{(s)}_{-2+\frac{r_{s}}{T}}w \equiv}{}
\times u^{(k_{2}-1)}_{-2+\frac{r_{k_{2}-1}}{T}}u^{(k_{2}+1)}_{-2+\frac{r_{k_{2}+1}}{T}}\cdots u^{(s)}_{-2+\frac{r_{s}}{T}}w \quad \bigl(\mod\ C^{T}_{3}(W)\bigr),
\end{align*}
which means
\[
u^{(1)}_{-2+\frac{r_{1}}{T}}u^{(2)}_{-2+\frac{r_{2}}{T}}\cdots u^{(s)}_{-2+\frac{r_{s}}{T}}w\in C^{T}_{3}(W).
\]

Suppose that no two of the original residues have sum at most $T-1$;
otherwise the preceding argument applies directly. Since $s=T+1$,
choose any ordering $u^{(1)},\dots,u^{(T+1)}$ and put
$a^{(1)}=u^{(1)}\in V^{q_1}$ with $q_1=r_1$.
At each stage $j\geq2$, if $q_{j-1}+r_j\leq T-1$, stop and set
$t=j-1$. The two required vectors are then $a^{(t)}\in V^{q_t}$
and $u^{(t+1)}\in V^{r_{t+1}}$; their residues have sum at most $T-1$.
If $q_{j-1}+r_j\geq T$, instead set
\smash{$a^{(j)}=a^{(j-1)}_{-1}u^{(j)}\in V^{q_j}$}, where
$q_j=q_{j-1}+r_j-T$. If this second case occurred at every stage
$2\leq j\leq T+1$, then
\[
q_{T+1}=\sum_{i=1}^{T+1}r_i-T^2
\leq(T+1)(T-1)-T^2=-1,
\]
contradicting $q_{T+1}\geq0$. Thus the process stops by stage
$T+1$. It cannot stop at $j=2$ under the present assumption, so
$2\leq t\leq T$. This proves the required existence, including
the boundary sums $2T-1$, $3T-2$, $3T-1$, and their analogues at later stages.
\begin{equation}\label{u(T+1)C3-2}
a^{(t)}\in V^{q_t},\qquad u^{(t+1)}\in V^{r_{t+1}},
\qquad 0\leq q_t+r_{t+1}\leq T-1.
\end{equation}

Next, for the nonnegative integer $t$ ($2\leq t\leq T$), assume the desired vectors are
\[
\begin{cases} \bigl(\bigl(\cdots\bigl(\bigl(u^{(1)}_{-1}u^{(2)}\bigr)_{-1}u^{(3)}\bigr)_{-1}\cdots\bigr)_{-1}u^{(t-1)}\bigr)_{-1}u^{(t)}\in V^{r} ,\\
u^{(t+1)}\in V^{r_{t+1}}\end{cases}
\]
for $u^{(i)}\in V^{r_{i}}$, $1\leq i\leq t+1$, $0\leq r_{i},r\leq T-1$ with $0\leq r+r_{t+1}\leq T-1$. For $0\leq q_{j}\leq T-1$, $1\leq j\leq t$, set
\[
a^{(j)}=\bigl(\bigl(\cdots\bigl(\bigl(u^{(1)}_{-1}u^{(2)}\bigr)_{-1}u^{(3)}\bigr)_{-1} \cdots\bigr)_{-1}u^{(j-1)}\bigr)_{-1}u^{(j)}\in V^{q_{j}},
\]
then it follows $a^{(j)}=a^{(j-1)}_{-1}u^{(j)}$. From the above proof of the existence of nonnegative integer~$t$ and vectors in (\ref{u(T+1)C3-2}), for any $j$ with ${2}\leq j\leq t$ we have $T\leq q_{j-1}+r_{j}\leq 2T-2$, thus $q_{j}=q_{j-1}+r_{j}-T$. By (\ref{twisted-ite-2}), there exist nonnegative integers $k_{1}$ and $l_{1}$ such that
\begin{gather}
 \bigl(a^{(t)}_{-1}u^{(t+1)}\bigr)_{-3+\frac{q_{t}+r_{t+1}}{T}}w\notag\\
 \qquad = \sum_{i=0}^{k_{1}}\sum_{j\in \mathbb{N}}\binom{-l_{1}-\frac{q_{t}}{T}}{i} \binom{i-1}{j}(-1)^{j}a^{(t)}_{-1+l_{1}+i-j+\frac{q_{t}}{T}}u^{(t+1)}_{-3-l_{1}-i+j+\frac{r_{t+1}}{T}}w.\label{u(T+1)C3-3}
\end{gather}
Since $0\leq q_{t}+r_{t+1}\leq T-1$, we have
\[
\bigl(a^{(t)}_{-1}u^{(t+1)}\bigr)_{-3+\frac{{q_{t}}+r_{t+1}}{T}}w\in C^{T}_{3}(W).
\]
If $i-1\geq 0$, then $j\leq i-1$. For each $i$, $j$, from (\ref{twisted-com}) we have
\begin{align*}
& a^{(t)}_{-1+l_{1}+i-j+\frac{q_{t}}{T}}u^{(t+1)}_{-3-l_{1}-i+j+\frac{r_{t+1}}{T}}w\\
& = u^{(t+1)}_{-3-l_{1}-i+j+\frac{r_{t+1}}{T}}a^{(t)}_{-1+l_{1}+i-j+\frac{q_{t}}{T}}w
 +\!\sum_{n\in \mathbb{N}}\!\binom{-1+l_{1}+i-j+\frac{q_{t}}{T}}{n}(u^{(t+1)}_{n}a^{(t)})_{-4-n+\frac{q_{t}+r_{t+1}}{T}}w.
\end{align*}
Notice that $-4-n\leq -3$ and $-3-l_{1}-i+j\leq -3$, we have
\begin{gather*}
\bigl(u^{(t+1)}_{n}a^{(t)}\bigr)_{-4-n+\frac{q_{t}+r_{t+1}}{T}}w\in C^{T}_{3}(W),\\
u^{(t+1)}_{-3-l_{1}-i+j+\frac{r_{t+1}}{T}}a^{(t)}_{-1+l_{1}+i-j+\frac{q_{t}}{T}}w\in C^{T}_{3}(W),
\end{gather*}
which implies
\[
a^{(t)}_{-1+l_{1}+i-j+\frac{q_{t}}{T}}u^{(t+1)}_{-3-l_{1}-i+j+\frac{r_{t+1}}{T}}w\in C^{T}_{3}(W).
\]

For $i=0$, if $j\geq l_{1}+2$, we have $-1+l_{1}-j\leq -3$, so that
\[
a^{(t)}_{-1+l_{1}-j+\frac{q_{t}}{T}}u^{(t+1)}_{-3-l_{1}+j+\frac{r_{t+1}}{T}}w\in C^{T}_{3}(W).
\]
If $j\leq l_{1}$, we have $-3-l_{1}+j\leq-3$. From (\ref{twisted-com}), we have
\begin{align*}
\begin{split}
& a^{(t)}_{-1+l_{1}-j+\frac{q_{t}}{T}}u^{(t+1)}_{-3-l_{1}+j+\frac{r_{t+1}}{T}}w\\
& \qquad = u^{(t+1)}_{-3-l_{1}+j+\frac{r_{t+1}}{T}}a^{(t)}_{-1+l_{1}-j+\frac{q_{t}}{T}}w
 +\sum_{n\in \mathbb{N}}\binom{-1+l_{1}-j+\frac{q_{t}}{T}}{n}(u^{(t+1)}_{n}a^{(t)})_{-4-n+\frac{q_{t}+r_{t+1}}{T}}w.
 \end{split}
\end{align*}
Thus we have
\[
a^{(t)}_{-1+l_{1}-j+\frac{q_{t}}{T}}u^{(t+1)}_{-3-l_{1}+j+\frac{r_{t+1}}{T}}w\in C^{T}_{3}(W).
\] So, for the only remaining case $i=0$ and $j=l_{1}+1$, the corresponding term $a^{(t)}_{-2+\frac{q_{t}}{T}}u^{(t+1)}_{-2+\frac{r_{t+1}}{T}}w$ must also lie in $C_{3}^{T}(W)$.

Next we set
\[
w'=u^{(t+1)}_{-2+\frac{r_{t+1}}{T}}w.
\] Since $a^{(t)}=a^{(t-1)}_{-1}u^{(t)}$ and $q_{t}=q_{t-1}+r_{t}-T$, then by (\ref{twisted-ite-2}) there exist nonnegative integers $k_{2}$ and $l_{2}$ such that
\begin{align*}
& \bigl(a^{(t-1)}_{-1}u^{(t)}\bigr)_{-2+\frac{q_{t}}{T}}w'
 \\
&\qquad = \sum_{i=0}^{k_{2}}\sum_{j\in \mathbb{N}}\binom{-l_{2}-\frac{q_{t-1}}{T}}{i}\binom{i-1}{j}(-1)^{j}a^{(t-1)}_{-1+l_{2}+i-j+\frac{q_{t-1}}{T}}u^{(t)}_{-2+\frac{q_{t}}{T}-l_{2}-i+j-\frac{q_{t-1}}{T}}w'\\
&\qquad = \sum_{i=0}^{k_{2}}\sum_{j\in \mathbb{N}}\binom{-l_{2}-\frac{q_{t-1}}{T}}{i}\binom{i-1}{j}(-1)^{j}a^{(t-1)}_{-1+l_{2}+i-j+\frac{q_{t-1}}{T}}u^{(t)}_{-3-l_{2}-i+j+\frac{r_{t}}{T}}w'.
\end{align*}
Refer to the proof of~(\ref{u(T+1)C3-3}), similarly for the only remaining case $i=0$ and $j=l_{2}+1$, the corresponding term
\[
a^{(t-1)}_{-2+\frac{q_{t-1}}{T}}u^{(t)}_{-2+\frac{r_{t}}{T}}w' =a^{(t-1)}_{-2+\frac{q_{t-1}}{T}}u^{(t)}_{-2+\frac{r_{t}}{T}}u^{(t+1)}_{-2+\frac{r_{t+1}}{T}}w
\]
must also lie in $C_{3}^{T}(W)$. Therefore, repeat the above process $t$ times, the only remaining term
\[
u^{(1)}_{-2+\frac{r_{1}}{T}}u^{(2)}_{-2+\frac{r_{2}}{T}}\cdots u^{(t)}_{-2+\frac{r_{t}}{T}}u^{(t+1)}_{-2+\frac{r_{t+1}}{T}}w
\]
must also lie in $C_{3}^{T}(W)$, which completes the whole proof.
\end{proof}

From Lemma \ref{u(T+1)C3}, we immediately have the following corollary.
\begin{Corollary}\label{U-C2-C3}
Let $V$ be a vertex algebra and let $W$ be a $g$-twisted $V$-module. Then
\[
u^{(1)}_{-2+\frac{r_{1}}{T}}u^{(2)}_{-2+\frac{r_{2}}{T}}\cdots u^{(s)}_{-2+\frac{r_{s}}{T}}C^{T}_{2}(W)\subset C^{T}_{3}(W)
\]
for $s\geq T$, $u^{(i)}\in V^{r_{i}}$, $w\in W$, $0\leq r_{i}\leq T-1$, $1\leq i\leq s$.
\end{Corollary}

The following result provides the relation between $u_{-k+\frac{p}{T}}C^{T}_{k}(W)$ and $C^{T}_{k+1}(W)$ for $k\geq 3$.

\begin{Lemma}\label{u-k-cktw}Let $V$ be a vertex algebra and let $W$ be a $g$-twisted $V$-module. Let $k\in \mathbb{Z}$ and $k\geq3$. Then
\[
u_{-k+\frac{p}{T}}C^{T}_{k}(W)\subset C^{T}_{k+1}(W)
\]
for $u\in V^{p}$, $0\leq p\leq T-1$.
\end{Lemma}

\begin{proof}Let $u\in V^{p}$, $0\leq p\leq T-1$, from Definition~\ref{dftmcn} we have $C^{T}_{n}(W)$ is linearly spanned by the vectors $v_{-n+\frac{q}{T}}w$, with $v\in V^{q}$, $0\leq q\leq T-1$. By (\ref{twisted-ite-2}), there exist nonnegative integers $t$ and $l$ such that
\begin{gather}
(u_{-1}v)_{-2k+1+\frac{p+q}{T}}w
 \nonumber\\
 \qquad{} =\sum_{i=0}^{t}\sum_{j\in \mathbb{N}}\binom{-l-\frac{p}{T}}{i}\binom{i-1}{j}(-1)^{j}u_{-1+l+i-j+\frac{p}{T}}v_{-2k+1-l-i+j+\frac{q}{T}}w.\label{lemma11-1}
\end{gather}

For $0\leq p+q \leq T-1$, $k\geq3$, we have $-2k+1\leq-k-1$ which implies
\[
(u_{-1}v)_{-2k+1+\frac{p+q}{T}}w\in C^{{T}}_{k+1}(W).
\]
 If $i-1\geq 0$, then $j\leq i-1$. For each $i$, $j$, from (\ref{twisted-com}) we have
\begin{align*}
& u_{-1+l+i-j+\frac{p}{T}}v_{-2k+1-l-i+j+\frac{q}{T}}w\\
&\qquad =v_{-2k+1-l-i+j+\frac{q}{T}}u_{-1+l+i-j+\frac{p}{T}}w+\sum_{r\in \mathbb{N}}\binom{i-j+l-1+\frac{p}{T}}{r}(u_{r}v)_{-2k-r+\frac{p+q}{T}}w.
\end{align*}
Notice that $-2k-r\leq -k-1$ and $-2k+1-l-i+j\leq -k-1$, we have
\begin{gather*}
(u_{r}v)_{-2k-r+\frac{p+q}{T}}w\in C^{T}_{k+1}(W),\\
v_{-2k+1-l-i+j+\frac{q}{T}}u_{-1+l+i-j+\frac{p}{T}}w\in C^{T}_{k+1}(W),
\end{gather*}
which implies
\[
u_{-1+l+i-j+\frac{p}{T}}v_{-2k+1-l-i+j+\frac{q}{T}}w\in C^{T}_{k+1}(W).
\]
For $i=0$, if $j\geq k+l$, we have $-1+l-j\leq-k-1$, so that
\[
u_{-1+l-j+\frac{p}{T}}v_{-2k+1-l+j+\frac{q}{T}}w\in C^{T}_{k+1}(W).
\]
If $j\leq k+l-2$, we have $-2k+1-l+j\leq-k-1$. From (\ref{twisted-com}), we get
\begin{align*}
& u_{-1+l-j+\frac{p}{T}}v_{-2k+1-l+j+\frac{q}{T}}w\\
&\qquad =v_{-2k+1-l+j+\frac{q}{T}}u_{-1+l-j+\frac{p}{T}}w+\sum_{r\in \mathbb{N}}\binom{-j+l-1+\frac{p}{T}}{r}(u_{r}v)_{-2k-r+\frac{p+q}{T}}w.
\end{align*}
Thus we have
\[
u_{-1+l-j+\frac{p}{T}}v_{-2k+1-l+j+\frac{q}{T}}w\in C^{T}_{k+1}(W).
\] Therefore, the only remaining case $i=0$ and $j=k+l-1$, the corresponding term $u_{-k+\frac{p}{T}}v_{-k+\frac{q}{T}}w$ must also lie in $C_{k+1}^{T}(W)$. This proves
$u_{-k+\frac{p}{T}}C^{T}_{k}(W)\subset C^{T}_{k+1}(W)$.

For $T\leq p+q\leq 2T-2$, we have
\[
(u_{-1}v)_{-2k+1+\frac{p+q}{T}}w=(u_{-1}v)_{-2k+2+\frac{p+q-T}{T}}w\in C^{T}_{k+1}(W),
\]
since $-2k+2\leq-k-1$ for $k\geq 3$. If $i-1\geq 0$, we have $j\leq i-1$. By (\ref{twisted-com}), we have
\begin{align*}
\begin{split}
& u_{-1+l+i-j+\frac{p}{T}}v_{-2k+1-l-i+j+\frac{q}{T}}w\\
&\qquad = v_{-2k+1-l-i+j+\frac{q}{T}}u_{-1+l+i-j+\frac{p}{T}}w+\sum_{r\in \mathbb{N}}\binom{i-j+l-1+\frac{p}{T}}{r}(u_{r}v)_{-2k-r+1+\frac{p+q-T}{T}}w.
\end{split}
\end{align*}
Since $-2k-r+1\leq-k-1$ and $-2k+1-l-i+j\leq-k-1$, we have
\begin{gather*}
(u_{r}v)_{-2k-r+1+\frac{p+q-T}{T}}w\in C^{T}_{k+1}(W),\\
v_{-2k+1-l-i+j+\frac{q}{T}}u_{-1+l+i-j+\frac{p}{T}}w\in C^{T}_{k+1}(W).
\end{gather*}

For $i=0$, if $j\geq k+l$, we have $-1+l-j\leq-k-1$, so that we obtain
\[
u_{-1+l-j+\frac{p}{T}}v_{-2k+1-l+j+\frac{q}{T}}w\in C^{T}_{k+1}(W).
\]
If $j\leq k+l-2$, we have $-2k+1-l+j\leq-k-1$. From (\ref{twisted-com}), we have
\begin{align*}
& u_{-1+l-j+\frac{p}{T}}v_{-2k+1-l+j+\frac{q}{T}}w\\
&\qquad =v_{-2k+1-l+j+\frac{q}{T}}u_{-1+l-j+\frac{p}{T}}w+\sum_{r\in \mathbb{N}}\binom{-j+l-1+\frac{p}{T}}{r}(u_{r}v)_{-2k-r+1+\frac{p+q-T}{T}}w.
\end{align*}
Thus we have
\[
u_{-1+l-j+\frac{p}{T}}v_{-2k+1-l+j+\frac{q}{T}}w\in C^{T}_{k+1}(W).
\]
 Similarly, the only remaining term $u_{-k+\frac{p}{T}}v_{-k+\frac{q}{T}}w$ in~(\ref{lemma11-1}) must also lie in $C^{T}_{k+1}(W)$.

To sum up, it proves $u_{-k+\frac{p}{T}}C^{T}_{k}(W)\subset C^{T}_{k+1}(W)$ for $k\in \mathbb{Z}$, $k\geq 3$, $u\in V^{p}$, $0\leq p\leq T-1$.
\end{proof}

\begin{Proposition}\label{p-etnw-cntw+2}Let $V$ be a vertex algebra, let $W$ be a $g$-twisted $V$-module and let $n\geq0$. Then
\[
u^{(1)}_{-k_{1}+\frac{r_{1}}{T}}u^{(2)}_{-k_{2}+\frac{r_{2}}{T}}\cdots u^{(s)}_{-k_{s}+\frac{r_{s}}{T}}w\in C^{T}_{n+3}(W)
\]
for $u^{(i)}\in V^{r_{i}}$, $0\leq r_{i}\leq T-1$, ${1}\leq i\leq s$, $w\in W$, $k_{1}, k_{2}, \dots, k_{s}\geq 2$, $s\geq(T+1)2^{n}$.\end{Proposition}

\begin{proof}
From (\ref{tmcn-relation-2}), we have
\[
u_{-i+\frac{p}{T}}C^{T}_{n+3}(W)\subset C^{T}_{n+3}(W)
\]
 for $u\in V^{p}$, $0\leq p\leq T-1$, $i\geq 1$, $n\geq0$, it suffices to prove the assertion for $s=(T+1)2^{n}$. Since
\[
u_{-k+\frac{p}{T}}=\frac{1}{\bigl(k-1-\frac{p}{T}\bigr)\bigl(k-2-\frac{p}{T}\bigr)\cdots \bigl(2-\frac{p}{T}\bigr)}\bigl(\D^{k-2}u\bigr)_{-2+\frac{p}{T}}
\]
for $k\geq3$, so we only need to prove the assertion for $k_{1}=k_{2}=\cdots=k_{s}=2$. We are going to use induction on $n$. If $n=0$, by Lemma~\ref{u(T+1)C3}, for $s\geq T+1$, we have
\[
u^{(1)}_{-2+\frac{r_{1}}{T}}u^{(2)}_{-2+\frac{r_{2}}{T}}\cdots u^{(s)}_{-2+\frac{r_{s}}{T}}w\in C^{T}_{3}(W).
\] Assume the assertion holds for $n=l$, some nonnegative integer. Set $s=(T+1)2^{l+1}$, $m=(T+1)2^{l}$. For $u^{(i)}\in V^{r_{i}}$, $0\leq r_{i}\leq T-1$, $1\leq i\leq s$, $w\in W$. By inductive hypothesis, we have
\[
u^{(m+1)}_{-2+\frac{r_{m+1}}{T}}u^{(m+2)}_{-2+\frac{r_{m+2}}{T}}\cdots u^{(s)}_{-2+\frac{r_{s}}{T}}w\in C^{T}_{l+3}(W).
\]
So that we have
\begin{gather}
u^{(1)}_{-2+\frac{r_{1}}{T}}u^{(2)}_{-2+\frac{r_{2}}{T}}\cdots u^{(s)}_{-2+\frac{r_{s}}{T}}w\in u^{(1)}_{-2+\frac{r_{1}}{T}}u^{(2)}_{-2+\frac{r_{2}}{T}}\cdots u^{(m)}_{-2+\frac{r_{m}}{T}}C^{T}_{l+3}(W).\label{p-en-cn+2-1}
\end{gather}
From Definition~\ref{dftmcn}, we know that $C^{T}_{l+3}(W)$ is linearly spanned by the vectors $a_{-l-3+\frac{r}{T}}w'$, with $a\in V^{r}$, $w'\in W$. Using (\ref{tmcn-relation-2}) and (\ref{tmcn-relation-3}), we have
\begin{gather}
 u^{(1)}_{-2+\frac{r_{1}}{T}}u^{(2)}_{-2+\frac{r_{2}}{T}}\cdots u^{(m)}_{-2+\frac{r_{m}}{T}}a_{-l-3+\frac{r}{T}}w'\nonumber\\
\qquad \equiv a_{-l-3+\frac{r}{T}}u^{(1)}_{-2+\frac{r_{1}}{T}}u^{(2)}_{-2+\frac{r_{2}}{T}}\cdots u^{(m)}_{-2+\frac{r_{m}}{T}}w'\quad \bigl(\mod\ C^{T}_{l+4}(W)\bigr).\label{p-en-cn+2-2}
\end{gather}
Furthermore, by inductive hypothesis, we have
\[
u^{(1)}_{-2+\frac{r_{1}}{T}}u^{(2)}_{-2+\frac{r_{2}}{T}}\cdots u^{(m)}_{-2+\frac{r_{m}}{T}}w'\in C^{T}_{l+3}(W),
\]
which together with Lemma \ref{u-k-cktw} and $l\geq 0$ gives
\[
a_{-l-3+\frac{r}{T}}u^{(1)}_{-2+\frac{r_{1}}{T}}u^{(2)}_{-2+\frac{r_{2}}{T}}\cdots u^{(m)}_{-2+\frac{r_{m}}{T}}w'\in a_{-l-3+\frac{r}{T}}C^{T}_{l+3}(W)\subset C^{T}_{l+4}(W).
\]
Then by (\ref{p-en-cn+2-2}) we have
\[
u^{(1)}_{-2+\frac{r_{1}}{T}}u^{(2)}_{-2+\frac{r_{2}}{T}}\cdots u^{(m)}_{-2+\frac{r_{m}}{T}}a_{-l-3+\frac{r}{T}}w'\in C^{T}_{l+4}(W),
\]
proving that
\[
u^{(1)}_{-2+\frac{r_{1}}{T}}u^{(2)}_{-2+\frac{r_{2}}{T}}\cdots u^{(m)}_{-2+\frac{r_{m}}{T}}C^{T}_{l+3}(W)\subset C^{T}_{l+4}(W).
\]
Therefore, by (\ref{p-en-cn+2-1}) we have
\[
u^{(1)}_{-2+\frac{r_{1}}{T}}u^{(2)}_{-2+\frac{r_{2}}{T}}\cdots u^{(s)}_{-2+\frac{r_{s}}{T}}w\in C^{T}_{l+4}(W).
\] This concludes the proof.
\end{proof}

The relationship between ${E^{T}_{n}(W)}$ and $C^{T}_{n}(W)$ in twisted case is described as follows.

\begin{Theorem}\label{relation-en-cn}Let $V$ be a vertex algebra, let $W$ be a $g$-twisted $V$-module and let $\mathcal{E}_{W}^{T}=\{E^{T}_{n}(W)\}_{n\in \mathbb{Z}}$ be the associated decreasing sequence. Then for any $n\geq2$
\begin{gather}
 C^{T}_{n}(W)\subset E^{T}_{(n-2)T+1}(W),\label{en-cn-relation-1}\\
 E^{T}_{m}(W)\subset C^{T}_{n}(W)\qquad \text{whenever} \ m\geq \max\bigl\{1,(n-2)T(T+1)2^{n-3}\bigr\}.\label{en-cn-relation-2}
\end{gather}
Furthermore,
\begin{gather}
\bigcap_{n\geq0}E^{T}_{n}(W)=\bigcap_{n\geq2}C^{T}_{n+2}(W).\label{en-cn-relation-3}
\end{gather}
\end{Theorem}

\begin{proof}
For $u\in V^{p}$, $0\leq p\leq T-1$, $n\geq 2$, $w\in W$, from Definitions \ref{dftmcn} and~\ref{dftmf}, we have
\[
u_{-n+\frac{p}{T}}w=u_{-1-(n-1)+\frac{p}{T}}w\in E^{T}_{(n-1)T-p}(W)\subset E^{T}_{(n-2)T+1}(W).
\]
This proves (\ref{en-cn-relation-1}). For $n=2$, from Lemma \ref{span-entw} we have $E^{T}_{1}(W)\subset C^{T}_{2}(W)$. For $n\geq3$, consider a generic spanning element of $E^{T}_{m}(W)$ with $m\geq 1$:
\[
X=u^{(1)}_{-1-k_{1}+\frac{r_{1}}{T}}u^{(2)}_{-1-k_{2}+\frac{r_{2}}{T}}\cdots u^{(s)}_{-1-k_{s}+\frac{r_{s}}{T}}w,
\]
where $s\geq1$, ${u^{(i)}}\in V^{r_{i}}$, $0\leq r_{i}\leq T-1$, $1\leq i\leq s$, $w\in W$, $k_{1},k_{2},\dots,k_{s}\geq1$, with
$k_{1}+k_{2}+\cdots+k_{s}-\frac{r_{1}+r_{2}+\cdots+r_{s}}{T}\geq\frac{m}{T}$. If there exists $j$ with $1\leq j\leq s$ such that $k_{j}\geq n-1$, by~(\ref{cntw-def}) we have
\[
u^{(j)}_{-1-k_{j}+\frac{r_{j}}{T}}W\subset C^{T}_{1+k_{j}}(W)\subset C^{T}_{n}(W),
\]
 then by (\ref{tmcn-relation-2}) we have $X\in C^{T}_{n}(W)$. If $s\geq(T+1)2^{n-3}$, then by Proposition \ref{p-etnw-cntw+2} we have $X\in C^{T}_{n}(W)$. Next we demonstrate at least one of $k_{i}\geq n-1$ and $s\geq(T+1)2^{n-3}$ holds. Assume for any $i$ with $1\leq i\leq s$, we have $k_{i}\leq n-2$ and $s\leq(T+1)2^{n-3}-1$, thus $k_{1}+k_{2}+\cdots+k_{s}\leq s(n-2)$. Since $0\leq r_{i}\leq T-1$, we have
\begin{align*}
(n-2)\frac{T(T+1)}{T}2^{n-3} &\leq \frac{m}{T}
 \leq k_{1}+k_{2}+\cdots+k_{s}-\frac{r_{1}+r_{2}+\cdots+r_{s}}{T}\\
&\leq s(n-2)
 \leq (n-2)\bigl((T+1)2^{n-3}-1\bigr).
\end{align*}
This leads to a contradiction. This proves (\ref{en-cn-relation-2}).

Combining (\ref{en-cn-relation-1}) and (\ref{en-cn-relation-2}), we have (\ref{en-cn-relation-3}).
\end{proof}

From Theorem~\ref{relation-en-cn}, we immediately have the following corollary.

\begin{Corollary}\label{e1=c2}Let $V$ be a vertex algebra and $W$ be a $g$-twisted $V$-module, we have
\[
E^{T}_{1}(W)=C^{T}_{2}(W).
\]
\end{Corollary}

\begin{proof}
From (\ref{en-cn-relation-1}), we have $C^{T}_{2}(W)\subset E^{T}_{1}(W)$. On the other hand, from (\ref{en-cn-relation-2}) we have $E^{T}_{1}(W)\subset C^{T}_{2}(W)$, which proves $E^{T}_{1}(W)=C^{T}_{2}(W)$.
\end{proof}

\section[Generating subspace of (1/T)N-graded twisted modules of vertex algebras]{Generating subspace of $\boldsymbol{\frac{1}{T}\mathbb{N}}$-graded twisted modules\\ of vertex algebras}\label{section6}

In this section, we shall use the twisted module $\mathrm{gr}_{\mathcal{E}}^{T}(W)$ of vertex Poisson algebra $\mathrm{gr}_{\mathcal{E}}^{T}(V)$ and the relation between $\{E_{n}^{T}(W)\}_{n\in\mathbb{Z}}$ and $\{C_{n}^{T}(W)\}_{n\in\mathbb{Z}_{\geq2}}$ to prove that for any twisted module~$W$ of a~vertex algebra $V$, $C_{2}$-cofiniteness implies $C_{n}$-cofiniteness for all $n\geq 2$. Next we employ $\mathrm{gr}_{\mathcal{E}}^{T}(W)$ to study generating subspaces of $\frac{1}{T}\mathbb{N}$-graded twisted modules of lower truncated $\mathbb{Z}$-graded vertex algebras.

The following gives the {bilinear action} of the differential algebra $A=\mathrm{gr}^{T}_{\mathcal{E}}(V)$ on $\mathrm{gr}^{T}_{\mathcal{E}}(W)$ {defined by~(\ref{tm-action})}.

\begin{Lemma}\label{tm-generation}Let $V$ be a vertex algebra and let $A=\mathrm{gr}^{T}_{\mathcal{E}}(V)$ be the vertex Poisson algebra obtained in Theorem~{\rm \ref{th-vpa}}, which is in particular a $T\mathbb{N}$-graded $($unital$)$ differential algebra. Let~$W$ be a~$g$\nobreakdash-twisted $V$-module, the associated graded vector space $\mathrm{gr}^{T}_{\mathcal{E}}(W)$ is equipped with the bilinear map
\begin{align*}
A\otimes\mathrm{gr}^{T}_{\mathcal{E}}(W)&\longrightarrow\mathrm{gr}^{T}_{\mathcal{E}}(W),\\
 \bigl(u+E^{T}_{(m+1)T}(V)\bigr)\otimes\bigl(w+E^{T}_{n+1}(W)\bigr)& \longmapsto\bigl(u+E^{T}_{(m+1)T}(V)\bigr)\cdot\bigl(w+E^{T}_{n+1}(W)\bigr)
\end{align*}
where
\begin{gather}
\bigl(u+E^{T}_{(m+1)T}(V)\bigr)\cdot\bigl(w+E^{T}_{n+1}(W)\bigr)=u_{-1+\frac{p}{T}}w+E^{T}_{mT+n-p+1}(W)\label{tm-action}
\end{gather}
for $u\in E^{T}_{mT}(V)\cap V^{p}$, $w\in E^{T}_{n}(W)$, with $m,n\in\mathbb{N}$, $0\leq p\leq T-1$. Moreover, $\mathrm{gr}^{T}_{\mathcal{E}}(W)$ is generated under the iterated action of $A$ via \eqref{tm-action} by its degree-zero subspace $E^{T}_{0}(W)/E^{T}_{1}(W)$, i.e., for $n\geq0$, the degree-$n$ subspace $E^{T}_{n}(W)/E^{T}_{n+1}(W)$ of $\mathrm{gr}^{T}_{\mathcal{E}}(W)$ is linearly spanned by the vectors
\begin{gather}
\partial^{k_{1}}\bigl(u^{(1)}+E^{T}_{T}(V)\bigr)\partial^{k_{2}}\bigl(u^{(2)}+E^{T}_{T}(V)\bigr) \cdots\partial^{k_{s}}\bigl(u^{(s)}+E^{T}_{T}(V)\bigr) \bigl(w+E^{T}_{1}(W)\bigr)\label{tm-generation1}
\end{gather}
for $1\leq s \leq n$, $u^{(i)}\in V^{r_{i}}$, $w\in W$, $k_{1}\geq k_{2}\geq\cdots\geq k_{s}\geq0$ with $k_{1}+k_{2}+\cdots+k_{s}-\frac{r_{1}+r_{2}+\cdots+r_{s}}{T}=\frac{n}{T}$, $0\leq r_{i}\leq T-1$, $1\leq i\leq s$.
\end{Lemma}

\begin{proof}Notice that Proposition {\ref{unw-entw}} guarantees that the operation given in
(\ref{tm-action}) is well defined. The case $n=0$ is trivial. For $n\geq1$, from (\ref{entw-entwnitr}) we have $E^{T}_{n}(W)$ is linearly spanned by the vectors $u_{-2-i+\frac{p}{T}}w$, where $u\in E_{0}^{T}(V)\cap V^{p}$, $w\in E_{n-(i+1)T+p}^{T}(W)$ for $0\leq i\leq\frac{n+p}{T}-1$, then we have
\begin{gather}
 u_{-2-i+\frac{p}{T}}w+E^{T}_{n+1}(W)\notag\\
\qquad =\frac{1}{\bigl(i+1-\frac{p}{T}\bigr)\bigl(i-\frac{p}{T}\bigr)\cdots \bigl(1-\frac{p}{T}\bigr)}\bigl(\D^{i+1}u\bigr)_{-1+\frac{p}{T}}w+E^{T}_{n+1}(W)\notag\\
\qquad=\frac{1}{\bigl(i+1-\frac{p}{T}\bigr) \bigl(i-\frac{p}{T}\bigr)\cdots\bigl(1-\frac{p}{T}\bigr)}\bigl(\D^{i+1}u+E^{T}_{(i+2)T}(V)\bigr) \bigl(w+E_{n-(i+1)T+p+1}^{T}(W)\bigr)\notag\\
\qquad=\frac{1}{\bigl(i+1-\frac{p}{T}\bigr)\bigl(i-\frac{p}{T}\bigr)\cdots \bigl(1-\frac{p}{T}\bigr)}\bigl(\partial^{i+1}\bigl(u+E^{T}_{T}(V)\bigr)\bigr) \bigl(w+E_{n-(i+1)T+p+1}^{T}(W)\bigr).\label{tm-generation-1}
\end{gather}
We are going to prove by induction on $n$. For $n=1$, from (\ref{entw-entwnitr}) we have $E^{T}_{1}(W)/E^{T}_{2}(W)$ is {linearly} spanned by $u_{-2-i+\frac{p}{T}}w+E^{T}_{2}(W)$ with $u\in E_{0}^{T}(V)\cap V^{p}$, $w\in E_{1-(i+1)T+p}^{T}(W)$ for $0\leq i\leq\frac{1+p}{T}-1$. Then it follows $i=0$ and $p=T-1$, we have $E^{T}_{1}(W)/E^{T}_{2}(W)$ is linearly spanned by vectors
\[
u_{-2+\frac{T-1}{T}}w+E^{T}_{2}(W)=\frac{1}{1-\frac{T-1}{T}}\bigl(\partial\bigl(u+E^{T}_{T}(V)\bigr)\bigr) \bigl(w+E^{T}_{1}(W)\bigr).
\]
We assume the assumption holds for $n-1$. From (\ref{tm-generation-1}), since $i\geq0$, $0\leq p\leq T-1$, we have $n-(i+1)T+p\leq n-1$, then it follows immediately from induction that $E^{T}_{n}(W)/E^{T}_{n+1}(W)$ is linearly spanned by those vectors in~(\ref{tm-generation1}). For $n\geq1$ and $1\leq j \leq s$, by~(\ref{entw-elements1}) we notice that $k_{j}$ and $r_{j}$ appear in pairs and $k_{j}\geq1$, $0\leq r_{j}\leq T-1$, we have $k_{j}-\frac{r_{j}}{T}\geq\frac{1}{T}$, which means $1\leq s\leq n$. This concludes the proof.
\end{proof}

The following result is from \cite{li-abel}.

\begin{Lemma}\label{vpa-span}Let $V$ be a vertex algebra and let $\mathrm{gr}^{T}_{\mathcal{E}}(V)$ be the vertex Poisson algebra obtained in Theorem~{\rm \ref{th-vpa}}. Then $\mathrm{gr}^{T}_{\mathcal{E}}(V)$ is linearly spanned by the vectors
\[
\partial^{k_{1}}\bigl(v^{(1)}+E^{T}_{T}(V)\bigr)\partial^{k_{2}} \bigl(v^{(2)}+E^{T}_{T}(V)\bigr)\cdots\partial^{k_{s}}\bigl(v^{(s)}+E^{T}_{T}(V)\bigr)
\]
for $s\geq 1$, $v^{(i)}\in V$, $k_{1}>k_{2}>\cdots >k_{s}\geq 0$. In particular, $\mathrm{gr}^{T}_{\mathcal{E}}(V)$ as a differential algebra is generated by the subspace $E^{T}_{0}(V)/E^{T}_{T}(V)(=V/C_{2}(V))$.
\end{Lemma}

From Lemma \ref{vpa-span} and the definition of $E^{T}_{n}(V)$, we have the following (see \cite{rrc2,bu-span,gn-rq,nt}).

\begin{Lemma}\label{c2-cn}Let $V$ be a vertex algebra. If $V$ is $C_{2}$-cofinite, then $V$ is $E_{n}$-cofinite and $C_{n+2}$-cofinite for any $n\geq 0$.
\end{Lemma}

Combining Theorem~\ref{relation-en-cn}, Corollary~\ref{e1=c2}, Lemmas~\ref{tm-generation} and~\ref{c2-cn}, we have the following.
\begin{Proposition}\label{p-c2cofninte}
Let $V$ be a vertex algebra and let $W$ be a $g$-twisted $V$-module. If $V$ and $W$ are $C_{2}$-cofinite,
then $W$ is $E_{n}$-cofinite and $C_{n+2}$-cofinite for all $n\geq 0$.
\end{Proposition}
\begin{proof}
We are going to use induction on $n$. Since $C_{2}(V)=E^{T}_{T}(V)$ and $C^{T}_{2}(W)=E^{T}_{1}(W)$, we have
\[
\dim V/E^{T}_{T}(V)= \dim V/C_{2}(V)<\infty.
\] From Lemma \ref{tm-generation}, we know that the degree-$n$ subspace $E^{T}_{n}(W)/E^{T}_{n+1}(W)$ is finite-dimensional for $n\geq 0$. For $n=0$, we have
$\dim E^{T}_{0}(W)/E^{T}_{1}(W)= \dim W/C^{T}_{2}(W)<\infty$.
 For $E^{T}_{0}(W)/E^{T}_{n+1}(W)$, we have
\[
E^{T}_{0}(W)/E^{T}_{{n}}(W)\cong \bigl(E^{T}_{0}(W)/E^{T}_{n{+1}}(W)\bigr)/\bigl(E^{T}_{n}(W)/E^{T}_{n+1}(W)\bigr).
\]
By induction hypothesis,
$\dim \bigl(E^{T}_{0}(W)/E^{T}_{n}(W)\bigr)<\infty$.
Then $\dim\bigl(E^{T}_{0}(W)/E^{T}_{n+1}(W)\bigr)<\infty$,
which proves $W$ is $E_{n}$-cofinite. By (\ref{en-cn-relation-2}),
$E^{T}_{m}(W)\subset C^{T}_{n}(W)$ when $m\geq (n-2)T(T+1)2^{n-3}$. Then
\[
\dim W/C^{T}_{n}(W)\leq\dim W/E^{T}_{m}(W)<\infty.
\]
This completes the proof.
\end{proof}

The following result generalizes the result of Lemma \ref{wt-entw}.

\begin{Proposition}\label{p-agtm-en-cn-relation}
Let $V=\coprod_{n\geq t}V_{(n)}$ be a lower truncated $\mathbb{Z}$-graded vertex algebra for some $t\in \mathbb{Z}$, and $W=\bigoplus_{n\in \frac{1}{T}\mathbb{N}}W_{(n)}$ be a $\frac{1}{T}\mathbb{N}$-graded $g$-twisted $V$-module. Then
\begin{gather}
C^{T}_{n}(W)\subset \coprod_{k\geq n+t-\frac{2T-1}{T}} W_{(k)}\qquad \mbox{for $n\geq 2$}.\label{p-agtm-en-cn-relation-1}
\end{gather}
Furthermore, for $n\geq 2$, we have
\begin{gather}
 \bigcap_{n\geq0}E^{T}_{n}(W)=\bigcap_{n\geq2}C^{T}_{n}(W)=0,\label{p-agtm-en-cn-relation-2}\\
 E^{T}_{m}(W)\subset \coprod_{k\geq n+t-\frac{2T-1}{T}}W_{(k)}\qquad \mbox{for $m\geq (n-2)T(T+1)2^{n-{3}}$.}\label{p-agtm-en-cn-relation-3}
\end{gather}
\end{Proposition}

\begin{proof}
For $n\geq2$ and for homogeneous vectors $u\in V^{p}$, $w\in W_{(m)}$, $m\in \frac{1}{T}\mathbb{N}$, we have
\[
\wt (u_{-n+\frac{p}{T}}w)=\wt u+\wt w +n-\frac{p}{T}-1\geq t+n-\frac{p}{T}-1\geq n+t-\frac{2T-1}{T}.
\]
Thus we prove (\ref{p-agtm-en-cn-relation-1}). Then we have $\bigcap_{n\geq2}C^{T}_{n}(W)=0$ which implies (\ref{p-agtm-en-cn-relation-2}). Combining Theorem~\ref{relation-en-cn} and (\ref{p-agtm-en-cn-relation-1}), we immediately have (\ref{p-agtm-en-cn-relation-3}).
\end{proof}

Next we consider a special case of $g$-twisted $V$-module whose associated decreasing sequence~$\mathcal{E}_{W}^{T}$ is trivial.

\begin{Lemma}\label{w=c2w}
Let $V$ be a vertex algebra and let $W$ be a $g$-twisted $V$-module. If $W=C_{2}^{T}(W)$, then
\[
E^{T}_{n}(W)=C^{T}_{n+2}(W)=W\qquad \mbox{for all $n\geq 0$.}
\]
\end{Lemma}

\begin{proof}
We are going to use induction on $n$. Since $W=C^{T}_{2}(W)$, by Corollary~\ref{e1=c2} we have
\[
W=C^{T}_{2}(W)=E^{T}_{1}(W)\subset W,
\] which implies $E^{T}_{1}(W)=W$. For some $k\geq 1$, assume
$E^{T}_{k}(W)=W$. For $v\in V^{p}$, $0\leq p\leq T-1$, then by Definitions~\ref{dftmcn} and~\ref{dftmf} we have
\[
v_{-2+\frac{p}{T}}W=v_{-2+\frac{p}{T}}E^{T}_{k}(W)\subset E^{T}_{k+T-p}(W)\subset E^{T}_{k+1}(W),
\]
which implies $W=C^{T}_{2}(W)\subset E^{T}_{k+1}(W)$, proving that $E^{T}_{k+1}(W)=W$. By induction, we have $E^{T}_{n}(W)=W$ for all $n\geq 0$. From (\ref{en-cn-relation-2}), we have $C^{T}_{n}(W)=W$ for all $n\geq 2$.
\end{proof}

\begin{Proposition}\label{can-not-graded}
Let $V=\bigsqcup_{n\geq t}V_{(n)}$ be a lower truncated vertex algebra for some $t\in \mathbb{Z}$, let $W$ be a nonzero $g$-twisted $V$-module such that $C^{T}_{2}(W)=W$. Then there does not exist a~lower truncated $\frac{1}{T}\mathbb{N}$-grading $W=\bigoplus_{n\in \frac{1}{T}\mathbb{N}}W_{(n)}$ with which $W$ becomes a $\frac{1}{T}\mathbb{N}$-graded $g$-twisted $V$-module.
\end{Proposition}

\begin{proof}
Suppose that $W$ is a $g$-twisted $V$-module such that $C^{T}_{2}(W)=W$, then by Lemma \ref{w=c2w} we have $W=C^{T}_{n+2}(W)$ for $n\geq0$, so that $W=\bigcap_{n\geq0}C^{T}_{n+2}(W)$. Furthermore, if $W$ carries a $\frac{1}{T}\mathbb{N}$-grading such that $W$ becomes a $\frac{1}{T}\mathbb{N}$-graded $g$-twisted $V$-module{, then} by (\ref{p-agtm-en-cn-relation-2}) we have $\bigcap_{n\geq0}C^{T}_{n+2}(W)=0$, so that $W=\bigcap_{n\geq0}C^{T}_{n+2}(W)=0$.
\end{proof}

The following proposition gives a spanning property of certain type for $\frac{1}{T}\mathbb{N}$-graded twisted modules of lower truncated {$\mathbb{Z}$-graded} vertex algebras.

\begin{Proposition}\label{prop:generation-of-W} Let $V=\coprod_{n\geq t}V_{(n)}$ be a $\mathbb{Z}$-graded vertex algebra for some $t\in\mathbb{Z}$ and let $U$ be a graded subspace of $V$ such that $V=U+C_{2}(V)$. Let \smash{$W=\bigoplus_{n\in \frac{1}{T}\mathbb{N}}W_{(n)}$} be a~$\frac{1}{T}\mathbb{N}$-graded $g$-twisted $V$-module and let $M$ be a graded subspace of $W$ such that $W=M+C^{T}_{2}(W)$. Then~$W$ is linearly spanned by the vectors
\begin{gather}
u^{(1)}_{-1-k_{1}+\frac{r_{1}}{T}}u^{(2)}_{-1-k_{2}+\frac{r_{2}}{T}}\cdots u^{(s)}_{-1-k_{s}+\frac{r_{s}}{T}}w\label{generation of W}
\end{gather}
for $u^{(i)}\in U\cap V^{r_{i}}$, $0\leq r_{i}\leq T-1$, $1\leq i\leq s$, $w\in M$, $k_{1}\geq k_{2}\geq\cdots\geq k_{s}\geq 0$, $s\geq 1$.
\end{Proposition}

\begin{proof}
Let $K^{T}(W)$ be a subspace of $W$, which is spanned by the vectors in (\ref{generation of W}). It is clear that $K^{T}(W)$ is a graded subspace. For $m\geq 0$, set
$K^{T}_{m}(W)=K^{T}(W)\cap E^{T}_{m}(W)$. For $u\in E^{T}_{lT}(V)\cap V^{p}$, $l\geq0$, $0\leq p\leq T-1$, set
\begin{align*}
&\partial^{(k)}(u+E^{T}_{(l+1)T}(V))= \frac{1}{\bigl(k-\frac{p}{T}\bigr) \bigl(k-1-\frac{p}{T}\bigr)\cdots\bigl(1-\frac{p}{T}\bigr)}\partial^{k}\bigl(u+E^{T}_{(l+1)T}(V)\bigr),\\
&\D^{(k)}u+E^{T}_{(l+k+1)T}(V)= \frac{1}{\bigl(k-\frac{p}{T}\bigr)\bigl(k-1-\frac{p}{T}\bigr)\cdots\bigl(1-\frac{p}{T}\bigr)} \bigl(\D^{k}u+E^{T}_{(l+k+1)T}(V)\bigr).
\end{align*}
Since $C^{T}_{2}(W)=E^{T}_{1}(W)$, then $W=M+C^{T}_{2}(W)=M+E^{T}_{1}(W)$. From~(\ref{tm-generation1}), for any $m\geq0$, $E^{T}_{m}(W)/E^{T}_{m+1}(W)$ is linearly spanned by the vectors
\[
\partial^{(k_{1})}\bigl(u^{(1)}+E^{T}_{T}(V)\bigr)\partial^{(k_{2})} \bigl(u^{(2)}+E^{T}_{T}(V)\bigr)\cdots\partial^{(k_{s})}\bigl(u^{(s)}+E^{T}_{T}(V)\bigr)\bigl(w+E^{T}_{1}(W)\bigr)
\]
for $s\geq 0$, $u^{(i)}\in U\cap V^{r_{i}}$, $w\in M$, $k_{1}\geq k_{2}\geq\cdots\geq k_{s}\geq0$ with $k_{1}+k_{2}+\cdots+k_{s}-\frac{r_{1}+r_{2}+\cdots+r_{s}}{T}=\frac{m}{T}$. By (\ref{grv-1}) and (\ref{grv-2}), we have
\begin{align*}
& \partial^{(k_{1})}\bigl(u^{(1)}+E^{T}_{T}(V)\bigr)\partial^{(k_{2})} \bigl(u^{(2)}+E^{T}_{T}(V)\bigr)\cdots\partial^{(k_{s})}\bigl(u^{(s)}+E^{T}_{T}(V)\bigr)\bigl(w+E^{T}_{1}(W)\bigr)\\
&\qquad = \bigl(\D^{(k_{1})}u^{(1)}+E^{T}_{k_{1}T+T}(V)\bigr) \bigl(\D^{(k_{2})}u^{(2)}+E^{T}_{k_{2}T+T}(V)\bigr)\cdots \\
&\qquad \quad {}\times \bigl(\D^{(k_{s})}u^{(s)}+E^{T}_{k_{s}T+T}(V)\bigr)\bigl(w+E^{T}_{1}(W)\bigr)\\
&\qquad = u^{(1)}_{-1-k_{1}+\frac{r_{1}}{T}}u^{(2)}_{-1-k_{2}+\frac{r_{2}}{T}}\cdots u^{(s)}_{-1-k_{s}+\frac{r_{s}}{T}}w+E^{T}_{m+1}(W).
\end{align*}
It follows that $E^{T}_{m}(W)=K^{T}_{m}(W)+E^{T}_{m+1}(W)$. Thus we have
\begin{align*}
W = E^{T}_{0}(W)
 = K^{T}_{0}(W)+K^{T}_{1}(W)+\cdots+K^{T}_{n}(W)+E^{T}_{n+1}(W) \subset K^{T}(W)+E^{T}_{n+1}(W)
\end{align*}
for any $n\geq0$. Since both $K^{T}(W)$ and $E^{T}_{n+1}(W)$ are graded {subspaces} of $W$ and by (\ref{p-agtm-en-cn-relation-2}) we have $\bigcap_{m\geq0}E^{T}_{m}(W)=0$, it forces
\[
W=K^{T}(W)=\sum_{i\geq0}K^{T}_{i}(W),
\]
proving the desired spanning property.
\end{proof}

\subsection*{Acknowledgements}

The authors sincerely thank the anonymous referees for their valuable comments and constructive suggestions, which have significantly helped improve this paper. J.-C.~Sun would like to {acknowledge} the support of National Natural Science Foundation of China (Nos.~12071276, 11931009 and 12226402) and Science and Technology Commission of Shanghai Municipality (No.~25ZR1401126).

\vspace{-1mm}

\pdfbookmark[1]{References}{ref}
\LastPageEnding

\end{document}